\documentclass[reqno]{amsart}
\usepackage{amsmath,amssymb,amsthm}
\usepackage{enumerate}
\usepackage{mathrsfs}
\usepackage[utf8]{inputenc}
\usepackage{xcolor}
\usepackage{mathtools}
\usepackage{hyperref}

\newtheorem{theorem}{Theorem}[section]
\newtheorem{corollary}[theorem]{Corollary}
\newtheorem{proposition}[theorem]{Proposition}
\newtheorem{lemma}[theorem]{Lemma}

\theoremstyle{definition}  
\newtheorem{remark}[theorem]{Remark}

\newtheorem*{question}{Question}

\newcommand{\C}{\mathbb{C}}
\newcommand{\B}{\mathbb{B}}
\newcommand{\R}{\mathbb{R}}
\newcommand{\Z}{\mathbb{Z}}
\newcommand{\inner}[2]{\langle#1,#2\rangle}
\newcommand{\inL}{\widetilde{\Delta}}

\DeclareMathOperator{\DO}{\mathcal{D}}

\DeclareMathOperator{\ran}{ran}
\newcommand{\ux}{u\chi_{\B_N}}

\newcommand{\Aw}{A^2_{\gamma}}
\newcommand{\Bw}{B_{\gamma}}
\newcommand{\Lw}{L^2_{\gamma}}
\newcommand{\Kw}{K^{\gamma}}
\newcommand{\dVw}{dV_{\gamma}}
\newcommand{\Hol}{\mathcal{H}(\B_N)}

\numberwithin{equation}{section}

\date{\today}

\begin{document}

\title[Brown--Halmos type theorems on the ball]{A generalization of the Brown--Halmos theorems for the unit ball}

\author{Trieu Le}
\address{Department of Mathematics and Statistics, University of Toledo, Toledo, OH 43606}
\email{Trieu.Le2@utoledo.edu}

\author{Akaki Tikaradze}
\address{Department of Mathematics and Statistics, University of Toledo, Toledo, OH 43606}
\email{Akaki.Tikaradze@utoledo.edu}

\begin{abstract}
In this paper we generalize the classical theorems of Brown and Halmos about algebraic properties of Toeplitz operators to Bergman spaces over the unit ball in several complex variables. A key result, which is of independent interest, is the characterization of summable functions $u$ on the unit ball whose Berezin transform can be written as a finite sum $\sum_{j}f_j\,\bar{g}_j$ with all $f_j, g_j$ being holomorphic.
In particular, we show that such a function must be pluriharmonic if it is sufficiently smooth and bounded. We also settle an open question about $\mathcal{M}$-harmonic functions. Our proofs employ techniques and results from function and operator theory as well as partial differential equations.
\end{abstract}



\maketitle

\section{Introduction and main results}

In their seminal work \cite{BH}, Brown and Halmos classified
all pairs of commuting Toeplitz operators on the Hardy space over the unit disc, as well as characterized all triples of Toeplitz operators $(T_f, T_g, T_h)$  such that $T_fT_g=T_h$. They showed that the product of two Toeplitz operators is zero if and only if one of them is zero. These theorems are commonly referred to as the Brown--Halmos theorems. Extending these results to the Bergman space setting and to Hilbert spaces of holomorphic functions on more general domains in several complex variables has been one of the central themes of research in the theory of Toeplitz operators in the last few decades.

On the Bergman space over the unit disc, the first results in the spirit of the Brown--Halmos theorems were obtained by Axler and {\v{C}}u\v{c}kovi\'{c} \cite{AxC} and Ahern and {\v{C}}u\v{c}kovi\'{c} \cite{AC}. It was shown in these papers that Brown--Halmos theorems hold true on the Bergman space for Toeplitz operators with bounded harmonic symbols. Subsequently, using his study of the range of the Berezin transform, Ahern \cite{A} improved the main result in \cite{AC}. Guo, Sun and Zheng \cite{GSZ} later studied finite rank semi-commutators and commutators of Toeplitz operators with harmonic symbols. It was shown that semi-commutators and commutators have finite rank if and only if they are actually zero. As a consequence, characterizations of the symbols were given. {\v{C}}u{\v{c}}kovi{\'c} \cite{CZ} obtained criteria for $T_{f}T_{g}-T_{h^n}$ to have finite rank, where $f,g$ and $h$ are bounded harmonic. More general results in this direction were investigated in \cite{CKL2008}. In a recent paper, Ding, Qin and Zheng \cite{DQZ} provided a more complete answer to the possible rank of $T_{f}T_{g}-T_{h}$ under the assumption that $f, g$ are bounded harmonic and $h$ is a $C^2$-functions and $(1-|z|^2)^2\Delta h$ is integrable. A complete characterization of these functions was then obtained. 
   
Researchers have also investigated Brown--Halmos theorems in the setting of several complex variables. A classification of pairs of commuting Toeplitz operators with pluriharmonic symbols on the unit ball was given by Zheng in \cite{Zh}. Subsequently, Choe and Koo \cite{CK} studied the zero product problem for Toeplitz operators on the unit ball with harmonic symbols having continuous extensions to part of the boundary. Finite sums of products of Toeplitz operators with pluriharmonic or $n$-harmonic symbols on the Bergman space over the polydisks were investigated in the papers by Choe et al. \cite{CLNZ2007, CKL2009}. The same problem on the Hardy space over the unit sphere was considered in \cite{CKLJMAA2011}. On the other hand, there has not {been} much progress in proving Brown--Halmos type results for Toeplitz operators with pluriharmonic symbols on the ball. It is our main goal to offer such results.
   
While it is not the focus of the current paper, we would like to mention that there is a vast literature on the study of Toeplitz operators with non-harmonic, non-pluriharmonic symbols. Researchers have investigated algebraic properties of Toeplitz operators whose symbols are radial, quasihomogeneous, or finite sums of quasihomogeneous functions in one and several variables on unweighted and standard weighted Bergman spaces. See, for example, \cite{App21,CL08,DZ16,GKV2003,GKV2004,JZD19,LR2008,Le08,LRZ2015,LSZ2006} and the references therein.
   
Throughout the paper, $N$ denotes a positive integer. We write $\B_N$ for the open unit ball in $\C^N$. We use $\Hol$ to the denote the algebra of all functions holomorphic on $\B_N$. For $\gamma>-1$, the weighted measure $dV_{\gamma}$ is defined by
\[dV_{\gamma}(z) = \frac{\Gamma(N+\gamma+1)}{N!\,\Gamma(\gamma+1)}(1-|z|^2)^{\gamma}dV(z),\] where $dV$ is the normalized Lebesgue measure on $\B_N$. Note that $dV_{\gamma}$ is a probability measure on $\B_N$. For $p>0$, the Bergman space $A^p_{\gamma}$ consists of all functions in $\Hol$ that are $p$-integrable with respect to $dV_{\gamma}$. The reader is referred to \cite[Chapter~2]{Zhu2005} for an excellent introduction to these Bergman spaces. We use $L^{p}_{\gamma}$ to denote the usual $L^{p}$-space with respect to $dV_{\gamma}$. Clearly, $A^{p}_{\gamma} = \Hol\cap L^{p}_{\gamma}$. Given $f\in L^{\infty}_{\gamma}$, one defines the Toeplitz operator $T_f:\Aw\longrightarrow \Aw$ by $T_{f}(h) = P_{\gamma}(fh)$ for all $h\in \Aw$. Here, $P_{\gamma}$ is the orthogonal projection from $\Lw$ onto $\Aw$. It is immediate that the operator $T_{f}$ is bounded and $\|T_{f}\|\leq\|f\|_{\infty}$. It is well known that $T_{f}$ can be expressed as an integral operator. In fact, we have
\begin{align*}
    T_{f}(h)(z) & = \int_{\B_N}\frac{f(w)h(w)}{(1-\inner{z}{w})^{N+1+\gamma}}\dVw(w),\quad z\in\B_N.
\end{align*}
For $f\in L^{1}_{\gamma}$, the above integral is also well defined for all bounded holomorphic functions $h$ on the ball. As a consequence, one may define Toeplitz operators for $L^{1}_{\gamma}$-symbols. It is well known that if $f\in L^1_{\gamma}$ is bounded on a set $\{z: r<|z|<1\}$ for some $0<r<1$, then $T_{f}$ extends to a bounded operator on $\Aw$. See the discussion on \cite[p.~204]{AC} for the one dimensional case. The general setting of several variables is similar.

We recall here some basic properties of Toeplitz operators. For bounded functions $\phi, \psi$ and complex numbers $a$ and $b$, we have
\[T_{a\phi+b\psi} = aT_{\phi}+bT_{\psi},\qquad T_{\phi}^{*} = T_{\bar{\phi}}.\]
It is also well known that if $\psi$ or $\bar{\phi}$ is holomorphic, then $T_{\phi}T_{\psi} = T_{\phi\psi}$. However, this property fails for general symbols, which is one of the reasons why the study of Toeplitz operators has attracted a great deal of attention. 

From now on, unless specified otherwise, we shall assume that $\gamma$ is a \textit{non-negative integer}. Our approach does not apply when $\gamma$ is not an integer. See Remark \ref{R:non_integer} for a more detailed discussion. The following theorem is the main result of the paper. It is a vast generalization of the aforementioned results and in a sense represents the best possible result one can hope for in the spirit of Brown--Halmos theorems for Toeplitz operators with pluriharmonic symbols. Recall that for two functions $x,y\in \Aw$, we use $x\otimes y$ to denote the operator \[(x\otimes y)(h) = \inner{h}{y}x,\quad h\in \Aw.\]

\begin{theorem}\label{BH}
Let $\phi_j, \psi_j$ be bounded pluriharmonic functions  for $1\leq j\leq n$ 
and $h$ be a bounded $C^{2N+2+2\gamma}$ function on $\B_{N}$. Let $x_{\ell}, y_{\ell}\in \Aw$ for $1\leq\ell\leq r$. 
Write  $\phi_j=f_j+\bar{g}_j, \psi_j=u_j+\bar{v}_j$ where $f_j, g_j, u_j$ and $v_j$ are holomorphic.
Then  
$$\sum_{j=1}^nT_{\phi_j}T_{\psi_j}=T_h+\sum_{\ell=1}^{r}x_{\ell}\otimes y_{\ell}$$ 
if and only if $h- \sum_{j=1}^{n}\bar{g}_ju_j$ is pluriharmonic and
\[\sum_{j=1}^{n} \phi_j\psi_j = h+(1-|z|^2)^{N+1+\gamma}\sum_{\ell=1}^{r}x_{\ell}\bar{y}_{\ell}.\]
\end{theorem}

As an immediate corollary, we have the following direct generalization of the Brown--Halmos theorems, which in particular settles the zero product problem for Toeplitz operators with pluriharmonic functions. The zero product problem for general symbols is a long standing open problem in the area of Toeplitz operators, which has resisted researchers' attempts even for the unit disc. Our result here in the single variable setting reduces to \cite[Theorem~7]{GSZ}.

\begin{corollary}\label{product}
Let $\phi, \psi$ be bounded pluriharmonic functions on $\B_N$. 

(a) If $T_{\phi}T_{\psi}=T_h$ on $\Aw$ for some bounded $C^{2N+2+2\gamma}$ function $h$, then $\bar{\phi}$ or $\psi$ is holomorphic and $\phi\psi=h$. 

(b) If $T_{\phi}T_{\psi}$ has a finite rank on $\Aw$, then $\phi$ or $\psi$ must be zero.
\end{corollary}

Another direct consequence of our main result is a strengthening of the aforementioned Zheng's theorem about commuting Toeplitz operators with pluriharmonic symbols. In the case of a single variable, we recover \cite[Theorem~6]{GSZ}.

\begin{corollary}\label{commutator}
Let $\phi, \psi$ be bounded pluriharmonic functions on $\B_N$. The commutator
$[T_{\phi}, T_{\psi}]$  has a finite rank on $\Aw$ if and only if both $\phi, \psi$ are holomorphic, or both are anti-holomorphic, or there are constants $c_1, c_2$, not both zero, such that $c_1\phi+c_2\psi$ is constant on $\B_N$.
\end{corollary}

The main tool for establishing our results is the Berezin transform. Recall that the Bergman space $\Aw$ is a reproducing kernel Hilbert space with kernel \[K^{\gamma}_z(w) = K^{\gamma}(w,z) = \frac{1}{(1-\inner{w}{z})^{N+1+\gamma}},\quad z,w\in\B_N.\]
Given a  function $u\in L^1_\gamma$, one defines the Berezin
transform of $u$ as follows 
\[\Bw(u)(z) = (1-|z|^2)^{N+1+\gamma}\int_{\B_N}\frac{u(\xi)}{|1-\inner{z}{\xi}|^{2(N+1+\gamma)}}dV_{\lambda}(\xi).\]
It is well known that
\[\Bw(u)(z) = \int_{\B_N}u\circ\varphi_{z}(\xi)\dVw(\xi),\]
where $\varphi_z$ is the automorphism of $\B_N$ that interchanges $0$ and $z$.
In the case $u\in L^2_\gamma$, we have
\[\Bw(u)(z) = \langle u k^{\gamma}_z, k^{\gamma}_z\rangle_{\Lw} = \langle T_{u}k^{\gamma}_z,k^{\gamma}_z\rangle,\]
where $k^{\gamma}_z=\Kw_{z}/\|\Kw_{z}\|$ is the normalized reproducing kernel. More generally, given a bounded operator $S:\Aw\to \Aw$, 
one defines similarly its Berezin transform 
$$\Bw(S)(z)=\langle S(k^{\gamma}_z), k^{\gamma}_z\rangle.$$
It is well known that the Berezin transform is an injective map. That is, if $\Bw(S_1)(z)=\Bw(S_2)(z)$ for all $z\in\B_N$, then $S_1=S_2$. The Berezin transform plays an important role in the theory of Toeplitz operators. In fact, it has been used as the main tool in the study of Brown--Halmos theorems for Toeplitz operators with pluriharmonic symbols in most of the references we have mentioned so far. 

It is clear that for $u\in L^{1}_{\gamma}$, the Berezin transform $\Bw(u)$ is real analytic on $\B_N$. As a result, we may expand $\Bw(u)$ as a series
\[\Bw(u)(z) = \sum_{\alpha, \beta}c_{\alpha,\beta}z^{\alpha}\bar{z}^{\beta}.\]
We say that $\Bw(u)$ has a \textit{finite rank} if the infinite matrix of coefficients $[c_{\alpha,\beta}]_{\alpha,\beta}$ has a finite rank. This happens if and only if there exist holomorphic functions $f_1,\ldots, f_n$ and $g_1,\ldots, g_n$ such that
\[\Bw(u)  = \sum_{j=1}^{n}f_j\,\bar{g}_j.\]

Besides being interesting on its own right, the following natural question is important in regards to algebraic properties of Toeplitz operators with pluriharmonic symbols.
\begin{question}
\label{Q:finite_rank_Berezin}
For which $u\in L^1_{\gamma}$ does $\Bw(u)$ have a finite rank?
\end{question}

For the unit disc on the complex plane, N.~V. Rao \cite{R} provided a full resolution of the above question in the unweighted case ($\gamma = 0$). Rao's result asserts that for an $L^1$-function $u$ on the unit disc, $B_{0}(u)$ has finite rank if and only if $u$ is harmonic except at a finite set of points. In particular, if $u$ is also assumed to be locally bounded, then it must be harmonic. In higher dimensions, the situation turns out to be more complicated and high dimensional phenomena do occur. In the theorem below, we describe $\Bw(u)$ whenever it is of finite rank. The proof of Theorem \ref{BH} relies heavily on this result. In addition, we answer an open question about $\mathcal{M}$-harmonic functions raised in \cite{CKLJMAA2011}. We recall here that $\mathcal{M}$-harmonic functions are those annihilated by the invariant Laplacian (see Section 2). It is well known that such functions are fixed points of the Berezin transforms.

\begin{theorem}\label{main_Berezin}
Suppose $u\in L^1_\gamma$ such that $B_{\gamma}(u)$ has a finite rank.
Then there exists a finite set $\Lambda\subset\overline{\B}_N$, a collection $\{P_{w}: w\in\Lambda\}$ of polynomials in $z$ and $\bar{z}$ of total degree at most $2N+1+2\gamma$, and a pluriharmonic function $h$ such that for $z\in\B_N$,
\[B_{\gamma}(u)(z) = h(z) + \sum_{w\in\Lambda}P_{w}\Big(\frac{z}{1-\inner{z}{w}}\Big).\]
Furthermore,
\begin{enumerate}[(a)]
    \item If $u$ also belongs to $C^{2N+2+2\gamma}(\B_N)$, then $\Lambda\subset \partial{\B}_{N}$, the unit sphere.
    \item If $u$ belongs to $L^{2N+2+2\gamma}_{\gamma}$, then  $\Lambda\subset\B_N$.
    \item If $B_{\gamma}(u) = f_1\bar{g}_1+\cdots+f_d\bar{g}_d$, where $f_{\ell}\in A^{2N+2+2\gamma}_{\gamma}$ and $g_{\ell}\in \Hol$ for all $\ell$, then $\Lambda\subset\B_N$. 
\end{enumerate}
As a consequence, if both (a) and (b), or both (a) and (c) hold, then $u$ is pluriharmonic.
\end{theorem}

\begin{remark}
As shown by Ahern and Rudin in \cite{AR} and can be verified directly, for $N\geq 3$, the function
\[u(z) = \frac{z_2\bar{z}_3}{|1-z_1|^2},\] which belongs to $L^1_{\gamma}$ whenever $\gamma\geq 0$, is $\mathcal{M}$-harmonic.
For such a function, we have $B_{\gamma}(u)=u$, so $B_{\gamma}(u)$ has a finite rank but $u$ is not pluriharmonic. 
In the case $N=2$, it can also be verified that 
\[u(z)=\frac{z_1\bar{z}_2}{(1-z_1)(1-\bar{z}_1)} -\frac{1}{2}\frac{\bar{z}_2^2{z}_2}{(1-z_1)(1-\bar{z}_1)^{2}}\]
is an $\mathcal{M}$-harmonic $L^1_{\gamma}$-function which is not pluriharmonic. It should be noted that Ahern and Rudin already showed that in the case of two complex variables, if $u=f\bar{g}$ (and $f,g$ are holomorphic) is $\mathcal{M}$-harmonic, then $u$ is actually pluriharmonic.
As a result, some type of regularity near the boundary
is required to conclude that $u$ is pluriharmonic, as in Theorem \ref{BH}. We would like to alert the reader that the existence 
of a smooth integrable function $u$ such that $B_{\gamma}(u)$ has a finite rank and $u$ is not
pluriharmonic is a high dimensional phenomenon. Indeed, it follows from the aforementioned result of Rao that if $u$ is \textit{locally bounded} on the unit disc (without any other assumption on regularity) so that $B_{0}(u)$ has a finite rank, then $u$ is harmonic.
\end{remark}

Our proof of Theorem \ref{main_Berezin} is influenced by Rao's idea to reformulate the finite rank property of the Berezin transform of $u$ in terms of a certain distribution associated to $u$ having a finite rank moment matrix, which allows the usage of a result due to Alexandrov and Rozenblum \cite{AlR}. In extending this approach to the case of the unit ball in  $\mathbb{C}^N$, significant complications do arise. We overcome these difficulties by establishing various identities for differential operators related to the invariant Laplacian and making use of a regularity result on integrable solutions of partial differential equations (see Section 2). 

Brown and Halmos proved that the zero operator is the only compact Toeplitz operator on the Hardy space over the unit disc. On Bergman spaces, there are many nontrivial compact Toeplitz operators. Indeed, whenever $f$ is a bounded function with a compact support contained in the unit ball, the operator $T_{f}$ is compact. On the other hand, the problem of determining nonzero finite rank Toeplitz operators was open for quite some time. In \cite{Lue}, Luecking settled this question in the negative by showing that whenever $\nu$ is a compactly supported finite measure on $\C$ for which the matrix of moments $[\int_{\C}z^{\ell}\bar{z}^{k}d\nu(z)]_{\ell,k}$ has finite rank, then $\nu$ is a linear combination of finitely many point masses. Luecking's theorem has been generalized to several complex variables \cite{Cho2009,RS2010} as well as to distributional symbols \cite{AlR}. We end this section by recalling the following result, which is crucial to our approach.

\begin{theorem}[Alexandrov-Rosenblum {\cite[Theorem~4.1]{AlR}}]\label{Luecking}
Let  $\mathcal{F}$ be a compactly supported distribution on
 $\mathbb{C}^N$. If the matrix  $\big[\mathcal{F}(\bar{z}^kz^{\ell})\big]_{k,\ell}$ has finite rank, then the support of $\mathcal{F}$ consists of finitely many points $\{w_1,\ldots, w_s\}$ and there are differential operators $L_1,\ldots, L_s$ such that $\mathcal{F} = \sum_{j=1}^{s}L_j(\delta_{w_j})$, where $\delta_w$ denotes the Dirac distribution at $w$.
\end{theorem}

\section{Some results on invariant Laplacian and radial derivative}

In this section we establish some results associated with certain
differential operators on the unit ball. Besides playing a crucial role in our  study of the Berezin transform, these identities are also interesting on their own right.

Ahern and {\v{C}}{u}{\v{c}}kovi{\'{c}} \cite{AC} and subsequently Ahern \cite{A}, Rao \cite{R} made of use the following property
of the kernel function of the Berezin transform (referred to as a ``marvelous identity" by Ahern)
\[\Delta_z\Big(\frac{(1-|z|^2)^2}{|1-z\bar{\xi}|^4}\Big)=\Delta_{\xi}\Big(\frac{(1-|\xi|^2)^2}{|1-z\bar{\xi}|^4}\Big).\]
The setting of several variables gets more complicated. We offer here several alternative identities which are important in our proofs. We shall make use of the following notation:
\[E_z = \sum_{j=1}^{N}z_j\frac{\partial}{\partial z_j},\qquad \bar{E}_z = \sum_{j=1}^{N}\bar{z}_j\frac{\partial}{\partial\bar{z}_j}, \qquad
\Delta_z = \sum_{j=1}^{N}\frac{\partial^2}{\partial z_j\partial\bar{z}_j}.\]
For any real number $s$, we write $|E_z+s|^2 = (E_z+s)(\bar{E}_z+s)$. We use $E_{\xi}, \bar{E}_{\xi}$, and $\Delta_{\xi}$ to denote the corresponding operators acting on the variable $\xi$. When the variable is not present, we write $E, \bar{E}$ and $\Delta$ for these operators.

To simplify the notation, for each integer $m\geq 1$, we define
$$\DO^{(m)}=(|E+m|^2-\Delta)\cdots (|E|^2-\Delta).$$
It is clear that $\DO^{(m)}$ is a differential operator of order $2m+2$ with polynomial coefficients. The following two lemmas provide important properties of $\DO^{(m)}$.

\begin{lemma}
\label{L:E_Delta_Identity}
For any integer $m\geq 1$ and $z,\xi\in\B_N$, we have
\begin{align}
\label{Eqn:E_Delta_Identity_A}
\DO^{(m)}_{\xi}\Big\{\frac{1}{|1-\inner{z}{\xi}|^{2}}\Big\}
    =  |E_{z}|^2\Big\{\frac{(m!)^2(1-|z|^2)^{m+1}}{|1-\inner{z}{\xi}|^{2(m+1)}}\Big\},
\end{align}
    and 
\begin{align}
\label{Eqn:E_Delta_Identity_B}
    \DO^{(m)}_{\xi}\Big\{\frac{1-|\xi|^2}{|1-\inner{z}{\xi}|^{2}}\Big\} = \big(|E_{z}|^2-\Delta_z\big)\Big\{\frac{(m!)^2(1-|z|^2)^{m+1}}{|1-\inner{z}{\xi}|^{2(m+1)}}\Big\}.
\end{align}
\end{lemma}

\begin{proof}
A direct calculation shows that
\begin{align}
\label{Eqn:E_Delta_a}
(|E_{\xi}|^2-\Delta_{\xi})\Big\{\frac{1}{|1-\inner{z}{\xi}|^2}\Big\} = \frac{|\inner{z}{\xi}|^2-|z|^2}{|1-\inner{z}{\xi}|^4}  = |E_z|^2\Big\{\frac{1-|z|^2}{|1-\inner{z}{\xi}|^2}\Big\},
\end{align}
\begin{align}
\label{Eqn:E_Delta_b}
 \Big(|E_{\xi}|^2-\Delta_{\xi}\Big)\Big\{\frac{1-|\xi|^2}{|1-\inner{z}{\xi}|^2}\Big\} & = \frac{(N-1)|1-\inner{z}{\xi}|^2+(1-|z|^2)(1-|\xi|^2)}{|1-\inner{z}{\xi}|^4}\notag\\
 & = \Big(|E_z|^2-\Delta_z\Big)\Big\{\frac{1-|z|^2}{|1-\inner{z}{\xi}|^2}\Big\},
\end{align}
and for any real number $s$,
\begin{align}
\label{Eqn:E_s_Delta}
\Big(|E_{\xi}+s|^2 - \Delta_{\xi}\Big)\Big\{\frac{1}{|1-\inner{z}{\xi}|^{2s}}\Big\} = \frac{s^2(1-|z|^2)}{|1-\inner{z}{\xi}|^{2(s+1)}}.
\end{align}
Applying $\big(|E_{\xi}+m|^2 - \Delta_{\xi}\big)\cdots(|E_{\xi}+1|^2-\Delta_{\xi})$ to \eqref{Eqn:E_Delta_a} and using \eqref{Eqn:E_s_Delta} repeatedly for $s=1,\ldots, m$ give \eqref{Eqn:E_Delta_Identity_A}. Finally, applying the same operator to \eqref{Eqn:E_Delta_b} and using \eqref{Eqn:E_s_Delta} give \eqref{Eqn:E_Delta_Identity_B}.
\end{proof}

We will use $\inL$, as usual, to denote the invariant Laplacian on $C^2(\B_N)$ which satisfies \[\inL = (1-|z|^2)(\Delta_{z}-|E_{z}|^2).\] Recall that $\inL$ can also be defined using the ordinary Laplacian and automorphisms of the unit ball. For more information on $\inL$ and its properties, see \cite[Chapter~4]{Ru} and \cite[Section~1.4]{Zhu2005}. However, the reader should be aware that the Laplacian defined there is actually four times our Laplacian. 

\begin{lemma}
\label{L:marvelous_identity}
For any integer $m\geq 1$ and $z,\xi\in\B_N$, we have
\begin{align}
\label{Eqn:marvelous_identity}
    (1-|\xi|^2)^{-m-1}p_{m}(\inL_{\xi})\left\{\frac{1-|\xi|^2}{|1-\inner{z}{\xi}|^2}\right\} & = (|E_z|^2-\Delta_z)\left\{\frac{(1-|z|^2)^{m+1}}{|1-\inner{z}{\xi}|^{2(m+1)}}\right\},
\end{align}
where \[p_m(t) = \frac{1}{(m!)^2}\prod_{j=0}^{m}\big(j(j-N)-t\big).\]
As a consequence,
\begin{align}
\label{Eqn:DO_identity}
(m!)^2(1-|\xi|^2)^{-m-1}p_{m}(\inL_{\xi}) = \DO^{(m)}_{\xi}.
\end{align}
\end{lemma}

\begin{proof}
Put $h(\xi)=1-|\xi|^2$. A direct but tedious calculation shows that for any positive integer $j\geq 1$,
\[j^2\,h^{j+1} = \big(j(j-N)-\inL\big)(h^j),\]
which implies
\[h^{m+1} = \frac{1}{(m!)^2}\prod_{j=1}^{m}\Big(j(j-N)-\inL\Big)(h).\]
Therefore,
\[-\inL h^{m+1} = \frac{1}{(m!)^2}(-\inL)\prod_{j=1}^{m}\Big(j(j-N)-\inL\Big)(h) = p_{m}(\inL)h.\]
Since both sides are radial functions that depend only on the modulus of the variable, for any $z,\xi\in\B_N$, we have
\[-(\inL h^{m+1})\circ\varphi_{\xi}(z) = (p_m(\inL)h)\circ\varphi_{z}(\xi).\]
On the other hand, the invariance of $\inL$ under the automorphisms of $\B_N$ gives
\[(\inL h^{m+1})\circ\varphi_{\xi}(z) = \inL_{z}(h^{m+1}\circ\varphi_{\xi}(z)) = \inL_z\big((1-|\varphi_{\xi}(z)|^2)^{m+1}\big),\]
and
\[ (p_m(\inL)h)\circ\varphi_{z}(\xi) = p_m(\inL_{\xi})(h\circ\varphi_{z}(\xi)) = p_m(\inL_{\xi})\big(1-|\varphi_{z}(\xi)|^2\big).\]
Consequently,
\begin{align}
\label{Eqn:p_m_inv_Lap}   
p_m(\inL_{\xi})\Big(1-|\varphi_{z}(\xi)|^2\Big) = -\inL_z\Big((1-|\varphi_{\xi}(z)|^2)^{m+1}\Big).
\end{align}
Since
\[1-|\varphi_z(\xi)|^2 = 1-|\varphi_{\xi}(z)|^2 = \frac{(1-|z|^2)(1-|\xi|^2)}{|1-\inner{z}{\xi}|^2}\]
and $-\inL_{z} = (1-|z|^2)(|E_z|^2-\Delta_z)$, the identity \eqref{Eqn:marvelous_identity} now follows from \eqref{Eqn:p_m_inv_Lap}.

From Lemma \ref{L:E_Delta_Identity} and equation \eqref{Eqn:marvelous_identity}, we conclude that the differential operators on both sides of  \eqref{Eqn:DO_identity} agree on functions of the form $\frac{1-|\xi|^2}{|1-\inner{z}{\xi}|^2}$ for all $z\in\B_N$. Taking partial derivatives in $z, \bar{z}$ and setting $z=0$, we see that the two operators agree on all polynomials of the form $(1-|\xi|^2)q(\xi,\bar{\xi})$, where $q$ is a polynomial. Since any smooth functions on $\B_N$ can be approximated by such polynomials (in the topology of uniform convergence on compact sets of all derivatives up to order $2m+2$), we obtain the required identity.
\end{proof}

{Let us now recall some basic definitions and background on the theory of distributions. There are many excellent resources but we shall use \cite[Chapters~II and III]{Hor03} as our main reference.

\begin{enumerate}[(1)]
\item Following \cite{Hor03}, we use $C^{\infty}_{0}(\R^n)$ to denote the space of $C^{\infty}$ functions on $\R^n$ having a compact support. For any subset $K\subset\R^n$, the space $C_{0}^{\infty}(K)$ consists of all functions in $C^{\infty}_{0}(\R^n)$ with support contained in $K$.

A distribution $\mathcal{F}$ in $\R^n$ of order at most $m$ is a linear functional on $C^{\infty}_{0}(\R^n)$ such that for each compact set $K\subset\R^n$, there exists a constant $C_{K}$ for which 
\[|\mathcal{F}(\phi)| \leq C_{K}\sum_{|\alpha|\leq m}\sup|\partial^{\alpha}\phi|\]
for all $\phi\in C^{\infty}_{0}(K)$. Here, $\alpha=(\alpha_1,\ldots,\alpha_n)\in\Z^{n}_{+}$, $|\alpha|=\alpha_1+\cdots+\alpha_n$, and $\partial^{\alpha} = \partial_{1}^{\alpha_1}\cdots\partial_n^{\alpha_n}$. The smallest value of $m$ for which the above condition holds is called the order of $\mathcal{F}$. 

\item To reduce the use of many parentheses in the following definitions, we shall write $\inner{\mathcal{F}}{\phi}$ to denote $\mathcal{F}(\phi)$. For any multiindex $\alpha$, the derivative $\partial^{\alpha}\mathcal{F}$ is a new distribution defined as
\[\inner{\partial^{\alpha}\mathcal{F}}{\phi} = (-1)^{|\alpha|}\inner{\mathcal{F}}{\partial^{\alpha}\phi},\qquad \phi\in C^{\infty}_{0}(\R^n).\]
For any smooth function $g$, the distribution $g\mathcal{F}$ is defined as
\[\inner{g\mathcal{F}}{\phi} = \inner{\mathcal{F}}{g\phi},\qquad \phi\in C_{0}^{\infty}(\R^n).\]
More generally, if $L$ is a differential operator of order $s$ with smooth coefficients written in the form 
\[L = \sum_{|\alpha|\leq s}a_{\alpha}(x)\partial^{\alpha},\] then
\[\inner{L(\mathcal{F})}{\phi} = \Big\langle\mathcal{F},\ \sum_{|\alpha|\leq s}(-1)^{|\alpha|}\partial^{\alpha}(a_{\alpha}\phi)\Big\rangle.\]
We define $L^{*}(\phi) = \sum_{|\alpha|\leq s}(-1)^{|\alpha|}\partial^{\alpha}(a_{\alpha}\phi)$ and call it the formal adjoint of $L$. The above formula can then be conveniently written as
\[\inner{L(\mathcal{F})}{\phi} = \inner{\mathcal{F}}{L^{*}(\phi)}.\]
Since the usual product rule holds true for products of smooth functions and distributions, we also have
\[\inner{L^{*}(\mathcal{F})}{\phi} = \inner{\mathcal{F}}{L(\phi)}.\]

\item The support of a distribution $\mathcal{F}$ is the set of points having no open neighborhood to which the restriction of $\mathcal{F}$ is zero. Section~2.3 in \cite{Hor03} discusses distributions with compact support. It is known that if $\mathcal{F}$ has a compact support, then its domain can be extended to the set of all functions that are smooth on an open neighborhood of the support of $\mathcal{F}$.

\cite[Theorem~2.3.4]{Hor03} asserts that if $\mathcal{F}$ is a distribution of order $k$ with support $\{a\}$, then there is a differential operator $L$ of order $k$ with constant coefficients such that $\mathcal{F} = L(\delta_a)$, where $\delta_a$ is the Dirac point mass distribution at $a$. It follows that if $\mathcal{F}$ has support $\{a_1,\ldots, a_s\}$, then 
\begin{align}
\label{Eqn:finitely_supported_dist}
\mathcal{F} = \sum_{j=1}^{s}L_j(\delta_{a_j})
\end{align} for some differential operators $L_1,\ldots, L_s$ with constant coefficients. The order of $\mathcal{F}$ equals to the maximum order of these operators.
\end{enumerate}
} 

We end the section with a result about singularities of $L^1$-solutions to PDEs. While we think this result may be known in the literature, due to the lack of an appropriate reference, we provide here a proof.
\begin{proposition}
\label{P:singularity_solutions_PDE}
Let $L$ be a differential operator of order $\mu$ with smooth coefficients on $\R^{n}$. Suppose $u\in L^1(\R^{n})$ having a compact support such that the distribution $L(u)$ is supported at finitely many points. Then the order of $L(u)$ is at most $\mu-1$.
\end{proposition}

\begin{remark}
In general, the distribution $L(u)$ may have order $\mu$. The point here is that if the support of $L(u)$ has only finitely many elements, then its order must be strictly smaller than $\mu$.
\end{remark}

\begin{proof}
Let $D_j = -i\,\partial_{j}$ and $D=(D_1,\ldots, D_n)$. Write the adjoint operator $L^{*}$ in the form
\begin{align}
    L^{*} = \sum_{|\alpha|\leq\mu}c_{\alpha}(x)D^{\alpha},
\end{align}
where each $c_{\alpha}$ is smooth. Since $u$ has a compact support, the Fourier transform of $L(u)$ can be computed by
\begin{align*}
   \widehat{L(u)}(\zeta) & = \int_{\R^n}u(x)L^{*}_{x}(e^{-i\inner{\zeta}{x}})dx\\
   & = \sum_{|\alpha|\leq\mu}\zeta^{\alpha}\int_{\R^n}u(x)c_{\alpha}(x)e^{-i\inner{\zeta}{x}}dx\\
   & = \sum_{|\alpha|\leq\mu}\zeta^{\alpha}\widehat{uc_{\alpha}}(\zeta).
\end{align*}
Because $c_{\alpha}$ is locally bounded, the function $uc_{\alpha}$ belongs to $L^1(\R^n)$. As a consequence, by the Riemann-Lebesgue Lemma, $\widehat{uc_{\alpha}}(\zeta)\rightarrow 0$ as $|\zeta|\rightarrow\infty$. It follows that for any $\zeta\in\R^n\backslash\{0\}$,
\begin{equation}
    \label{Eqn:asymptotic_Fourier_Lu}
    \lim_{t\rightarrow\infty}\frac{\widehat{L(u)}(t\zeta)}{t^{\mu}} = 0.
\end{equation}
Now let $\{a_1,\ldots,a_s\}$ be the support of $L(u)$. Then by \eqref{Eqn:finitely_supported_dist}, there are differential operators $L_1, \ldots, L_s$ with constant coefficients such that
\[L(u) = \sum_{j=1}^{s}L_j(\delta_{a_j}),\] where $\delta_a$ denotes the Dirac distribution at $a$. Since $L(u)$ has order at most $\mu$, each $L_j$ has order at most $\mu$ as well. As a result, there are homogeneous polynomials $p_j$ of degree $\mu$ such that $L_j = p_j(D) + \text{lower order derivatives}$. We then have
\begin{align*}
    \widehat{L(u)}(\zeta) & = \sum_{j=1}^{s}\widehat{L_j(\delta_{a_j})}(\zeta)\\
    & = \sum_{j=1}^{s}e^{-i\inner{\zeta}{a_j}}\cdot\big(p_j(\zeta) + \text{lower order terms in $\zeta$}\big).
\end{align*}
It now follows from \eqref{Eqn:asymptotic_Fourier_Lu} that
\begin{align}
\label{Eqn:limit_Fourier_transform}
\lim_{t\to\infty}\sum_{j=1}^{s}e^{-it\inner{\zeta}{a_j}}p_j(\zeta) = 0
\end{align}
for all $\zeta\in\R^n\backslash\{0\}$.

\underline{Claim}: for all $\zeta\in\R^n$ such that $\inner{\zeta}{a_j}\neq\inner{\zeta}{a_k}$ for all $j\neq k$, \eqref{Eqn:limit_Fourier_transform} forces $p_j(\zeta)=0$.

Since the set of all $\zeta$ in the claim  is dense in $\R^n$, we conclude that $p_j=0$ and hence  $L_j$ is of order at most $\mu-1$ for all $j$. Consequently, $L(u)$ has order at most $\mu-1$.

\textit{Proof of the claim.} We believe that the claim should be well known but we sketch here a proof.  To simplify the notation, put $\lambda_j=-\inner{\zeta}{a_j}$ and $b_j=p_j(\zeta)$. Note that the values $\lambda_1,\ldots,\lambda_s$ are pairwise distinct so there exists a real number $c$ such that $e^{i\lambda_1c},\ldots, e^{i\lambda_s c}$ are pairwise distinct.  Define $f(t)=\sum_{j=1}^{s}b_je^{i\lambda_jt}$ for $t\in\R$. Then \eqref{Eqn:limit_Fourier_transform} gives $\lim_{t\to\infty}f(t) = 0$ and hence, $\lim_{t\to\infty}f(t+\ell c)=0$ for all $0\leq \ell\leq s$. Note that
\[f(t+\ell c) = \sum_{j=1}^{s}(e^{i\lambda_j c})^{\ell}b_je^{i\lambda_jt}\]
so each $b_je^{i\lambda_jt}$ can be expressed as a linear combination of $f(t), f(t+c),\ldots, f(t+(s-1)c)$ via the Vandermonde determinant. It then follows that for each $1\leq j\leq s$, we have $\lim_{t\to\infty}b_je^{i\lambda_jt}=0$, which implies $b_j=0$.
\end{proof}

\section{Finite rank Berezin transform}
The goal of this section is to study finite rank Berezin transform $\Bw(u)$, which can be written in the form $\Bw(u)=\sum_{j=1}^{n}f_j\bar{g}_j$ for holomorphic functions $f_j$ and $g_j$. 
For a function $u\in L^1_\gamma$, we use $\ux$ to denote the corresponding distribution on $\C^N$ defined as
\[\phi\mapsto\int_{\B_N}u(\xi)\phi(\xi)\,\dVw(\xi),\quad \phi\in C^{\infty}(\C^N).\]
To apply the theory of distributions, we shall identity $\C^N$ with $\R^{2N}$.

The proof of Theorem \ref{main_Berezin} is divided into several steps. We now prove the first part and statement (a) in the theorem.

\begin{proposition}
\label{T:finite_rank_Bu}
Let $u\in L^1_{\gamma}$ such that $\Bw(u)$ has finite rank.
Then there exists a finite set $\Lambda\subset\overline{\B}_N$, a collection $\{P_{w}: w\in\Lambda\}$ of polynomials in $z$ and $\bar{z}$ of total degree at most $2N+1+2\gamma$, and a pluriharmonic function $h$ such that for $z\in\B_N$,
\[\Bw(u)(z) = h(z) + \sum_{w\in\Lambda}P_{w}\Big(\frac{z}{1-\inner{z}{w}}\Big).\]
If, furthermore, $u$ also belongs to $C^{2N+2+2\gamma}(\B_N)$, then $\Lambda\subset\partial\B_{N}$, the unit sphere.
\end{proposition}

\begin{proof}
Define
{\[\DO = \frac{1}{[(N+\gamma)!]^2}\DO^{(N+\gamma)} = \frac{1}{[(N+\gamma)!]^2}\Big(\big|E+(N+\gamma)\big|^2-\Delta\Big)\cdots \big(|E|^2-\Delta\big).\]}
Recall the formula for the Berezin transform
\[\Bw(u)(z) = \int_{\B_N}u(\xi)\frac{(1-|z|^2)^{N+1+\gamma}}{|1-\inner{z}{\xi}|^{2(N+1+\gamma)}}\dVw(\xi).\]
Applying $|E_z|^2$ to both sides and using Lemma \ref{L:E_Delta_Identity}, we conclude that
\begin{align}
|E_z|^2\Big(\Bw(u)(z)\Big) 
& = \int_{\B_N}u(\xi)\DO_{\xi}\Big\{\frac{1}{|1-\inner{z}{\xi}|^{2}}\Big\}\,\dVw(\xi)\label{Eqn:dist_E_Bu_a}\\
& = \int_{\B_N}u(\xi)\sum_{k, \ell}\binom{|k|}{k}\binom{|\ell|}{\ell}z^k\bar{z}^{\ell}\DO_{\xi}(\bar{\xi}^k{\xi}^{\ell})\,\dVw(\xi)\notag\\
& = \sum_{k,\ell}\Big\{\binom{|k|}{k}\binom{|\ell|}{\ell}\int_{\B_N}u(\xi)\DO_{\xi}(\bar{\xi}^{k}\xi^{\ell})\,\dVw(\xi)\Big\}z^k\bar{z}^{\ell}.\label{Eqn:dist_E_Bu_b}
\end{align}
Since $\Bw(u)$ is real analytic, we may write
\[\Bw(u)(z) = \sum_{k,\ell}a_{k,\ell}z^k\bar{z}^{\ell},\]
which implies
\[|E_z|^2(\Bw(u)(z)) = \sum_{k,\ell}|k||\ell|a_{k,\ell}z^k\bar{z}^{\ell}.\]
It follows that $\Bw(u)$ has a finite rank if and only if $|E_z|^2(\Bw(u))$ has a finite rank. Using \eqref{Eqn:dist_E_Bu_b}, we conclude that
\begin{align*}
   \Bw(u) \text{ has finite rank } 
   & \Longleftrightarrow  \Big[\binom{|k|}{k}\binom{|\ell|}{\ell}\int_{\B_N}u(\xi)\DO_{\xi}(\bar{\xi}^{k}\xi^{\ell})\,\dVw(\xi)\Big]_{k,\ell} \text{ has finite rank}\\
   & \Longleftrightarrow
   \Big[\int_{\B_N}u(\xi)\DO_{\xi}(\bar{\xi}^{k}\xi^{\ell})\,\dVw(\xi)\Big]_{k,\ell} \text{ has finite rank}.
\end{align*}
Consider the distribution $\mathcal{F}=\DO^{*}(u\chi_{\B_N})$, which is given by the formula
\[\mathcal{F}(\phi) =  \int_{\B_N}u(\xi)\DO(\phi)\dVw(\xi),\] for any $\phi$ that is smooth on an open neighborhood of $\overline{\B}_N$. Note that the support of $\mathcal{F}$ is contained in $\overline{\B}_N$. Since $\Bw(u)$ is assumed to have finite rank by the hypothesis, the matrix
\[\Big[\mathcal{F}(\bar{\xi}^{k}\xi^{\ell})\Big]_{k,\ell} = \Big[\int_{\B_N}u(\xi)\DO_{\xi}(\bar{\xi}^{k}\xi^{\ell})\,\dVw(\xi)\Big]_{k,\ell}\] has finite rank.
Applying Theorem \ref{Luecking}, we conclude that $\mathcal{F}$ has a finite support. 

Let $\{w_1,\ldots,w_s\}\subset\overline{\B}_N$ be the support of $\mathcal{F}$. By Proposition \ref{P:singularity_solutions_PDE}, $\mathcal{F}=\DO^{*}(\ux)$ has order at most $2N+1+2\gamma$ since $\DO^{*}$ is a differential operator of order $2N+2+2\gamma$. From formula \eqref{Eqn:dist_E_Bu_a}, we see that there are complex constants $a_{j,\alpha,\beta}$ for $1\leq j\leq s$ and $|\alpha|+|\beta|\leq 2N+1+\gamma$ such that
\begin{align}
\label{Eqn:formula_Ez_2_Bu}
|E_z|^2\Big(\Bw(u)(z)\Big)
& = \DO^{*}(\ux)\Big\{\frac{1}{|1-\inner{z}{\cdot}|^2}\Big\}\notag\\
& =\sum_{\substack{1\leq j\leq s\\ |\alpha|+|\beta|\leq 2N+1+\gamma }}a_{j,\alpha,\beta}\frac{z^{\alpha}\bar{z}^{\beta}}{(1-\inner{z}{w_j})^{1+|\alpha|}(1-\inner{w_j}{z})^{1+|\beta|}}.
\end{align}
A direct calculation shows that for any $w\in\overline{\B}_{N}$ and $|\alpha|, |\beta|\geq 1$, 
\begin{align*}
\frac{1}{1-\inner{z}{w}} - 1 & = \frac{\inner{z}{w}}{1-\inner{z}{w}} = E_z\Big\{\log\frac{1}{1-\inner{z}{w}}\Big\},\\
\frac{z^{\alpha}}{(1-\inner{z}{w})^{1+|\alpha|}} & = E_z\Big\{\frac{1}{|\alpha|}\frac{z^{\alpha}}{(1-\inner{z}{w})^{|\alpha|}}\Big\}.
\end{align*}
Thus, for such $\alpha, \beta$, the functions  $$\frac{z^{\alpha}\bar{z}^{\beta}}{(1-\inner{z}{w})^{1+|\alpha|}\,(1-\inner{w}{z})^{1+|\beta|}},$$  
$$\Big(\frac{1}{1-\inner{z}{w}}-1\Big)\frac{\bar{z}^{\beta}}{(1-\inner{w}{z})^{1+|\beta|}},\quad \frac{{z}^{\alpha}}{(1-\inner{z}{w})^{1+|\alpha|}}\Big(\frac{1}{1-\inner{w}{z}}-1\Big)$$
belong to the range of $|E_z|^2$. Hence, the identity \eqref{Eqn:formula_Ez_2_Bu} implies that the pluriharmonic function
\begin{align*}
& \sum_{\substack{1\leq j \leq s\\ |\alpha|\geq 1}}
a_{j,\alpha, 0}\frac{z^{\alpha}}{(1-\inner{z}{w_j})^{1+|\alpha|}}+\sum_{\substack{1\leq j\leq s\\ |\beta|\geq 1}}
a_{j,0, \beta}\frac{\bar{z}^{\beta}}{(1-\inner{w_j}{z})^{1+|\beta|}}\\
&\qquad + \sum_{0\leq j\leq s}a_{j,0,0}\Big(\frac{1}{1-\inner{z}{w_j}}+\frac{1}{1-\inner{w_j}{z}}-1\Big)
\end{align*}
is the image, under $|E_z|^2$, of a real analytic function. Using power series, we see that zero is the only pluriharmonic function belonging to the range of $|E_z|^2$. It then follows that for all $j$, we have $a_{j,\alpha,\beta}=0$ whenever $|\alpha|=0$ or $|\beta|=0$. Consequently,
\begin{align*}
    |E_z|^2\Big(\Bw(u)\Big) & = \sum_{\substack{1\leq j\leq s\\ |\alpha|\geq 1, |\beta|\geq 1}}a_{j,\alpha,\beta}\frac{z^{\alpha}\bar{z}^{\beta}}{(1-\inner{z}{w_j})^{1+|\alpha|}(1-\inner{w_j}{z})^{1+|\beta|}}\\
    & = |E_z|^2\Big\{\sum_{\substack{1\leq j\leq s\\ |\alpha|\geq 1, |\beta|\geq 1}}\frac{a_{j,\alpha,\beta}}{|\alpha|\,|\beta|}\frac{z^{\alpha}\bar{z}^{\beta}}{(1-\inner{z}{w_j})^{|\alpha|}(1-\inner{w_j}{z})^{|\beta|}}\Big\},
\end{align*}
which gives
\[\Bw(u)(z) = h(z)+\sum_{\substack{1\leq j\leq s\\ |\alpha|\geq 1, |\beta|\geq 1}}\frac{a_{j,\alpha,\beta}}{|\alpha|\,|\beta|}\frac{z^{\alpha}\bar{z}^{\beta}}{(1-\inner{z}{w_j})^{|\alpha|}(1-\inner{w_j}{z})^{|\beta|}},\] for some pluriharmonic function $h$ on $\B_N$. Defining 
\[P_{w_j}(z) = \sum_{\substack{|\alpha|\geq 1, |\beta|\geq 1\\ |\alpha|+|\beta|\leq 2N+1+\gamma}}\frac{a_{j,\alpha,\beta}}{|\alpha|\,|\beta|}z^{\alpha}\bar{z}^{\beta},\]
we obtain the required representation for $\Bw(u)$.

If $u$ belongs to $C^{2N+2+2\gamma}(\B_N)$, then the support of $\DO^{*}(\ux)$, being a finite set of points, must be contained in the unit sphere. As a result, $|w_j|=1$ for all $1\leq j \leq s$.
\end{proof}

From the proof of Proposition \ref{T:finite_rank_Bu} we have the following result which might be of independent interest.

\begin{corollary}
Let  $u\in L^1(\B_N,dV)$ be of the form  $u=\sum_{j=1}^{n}f_j\,\bar{g}_j$  with holomorphic $f_j, g_j$.  If  $u$ is an eigenfunction of the invariant Laplacian with eigenvalue $\lambda$, then
 $\lambda=j(j-N)$  for some $j\in \lbrace 0, 1,\ldots, N\rbrace.$ 

\end{corollary}
\begin{proof}
It is well known \cite[Theorem~4.2.4]{Ru} that eigenfunctions of  $\inL$ are also eigenfunctions of $B_0$ (the unweighted Berezin transform). Therefore, $B_0(u)$  has a finite rank. From the proof of Proposition \ref{T:finite_rank_Bu} as above (in the case $\gamma = 0$), we have that  $\DO^{*}(u)=0$ on $\B_N$, where
\begin{align*}
\DO & = \frac{1}{(N!)^2}\Big(\big|E+N\big|^2-\Delta\Big)\cdots \big(|E|^2-\Delta\big) = \frac{1}{(N!)^2}(1-|\xi|^2)^{-(N+1)}\prod_{j=0}^N\Big(j(j-N)-\inL\Big).
\end{align*}
The second equality comes from Lemma \ref{Eqn:DO_identity}. Since
\[(E+j)^{*} = -E + (j-N),\quad (\bar{E}+j)^{*} = -\bar{E} + (s-N) \text{ and } (\Delta)^{*} = \Delta\] for all $0\leq j\leq N$, we conclude that $\DO^{*}  = \DO$. The desired result then follows.
\end{proof}

Applying Proposition \ref{T:finite_rank_Bu} to the case where $u$ belongs to $L^{2N+2+2\gamma}_{\gamma}$, we now prove statement (b) in Theorem \ref{main_Berezin}.

\begin{proposition}
\label{T:u_L2N}
Suppose $u\in L^{2N+2+2\gamma}_{\gamma}$ and $\Bw(u)$ has finite rank. Then there exist finitely many points $w_1,\ldots,w_s\in \B_{N}$, polynomials $Q_1,\ldots, Q_s$ in $\C[z,\bar{z}]$ with total degrees at most $2N+1+2\gamma$, and a pluriharmonic function $h$ such that
\[\Bw(u)(z) = h(z) + \sum_{j=1}^{s}Q_{j}\circ\varphi_{w_j}(z).\]
\end{proposition}

\begin{proof}
We know that there exist holomorphic functions $h_1, h_2$ on $\B_N$ and finitely many points $w_1,\ldots,w_s\in\overline{\B}_N$ such that
\begin{align}
\label{Eqn:B_u_formula}
\Bw(u)(z) & = h_1(z)+\bar{h}_2(z) + \sum_{\substack{1\leq j\leq s\\
1\leq |\beta|\leq 2N+2\gamma}}Q_{j,\beta}\Big(\frac{z}{(1-\inner{z}{w_j})}\Big)\,\cdot\frac{\bar{z}^{\beta}}{(1-\inner{w_j}{z})^{|\beta|}},
\end{align}
where each $Q_{j,\beta}$ is a holomorphic polynomial of degree at most $2N+1+2\gamma-|\beta|$ with $Q_{j,\beta}(0)=0$. We prove first that $Q_{j,\beta}=0$ whenever $|w_j|=1$.

Complexifying \eqref{Eqn:B_u_formula} gives
\begin{align*}
    \Bw(u)(z,\zeta)  - \bar{h}_2(\zeta) & = h_1(z)+ \sum_{\substack{1\leq j\leq s\\
1\leq |\beta|\leq 2N+2\gamma}}Q_{j,\beta}\Big(\frac{z}{(1-\inner{z}{w_j})}\Big)\,\cdot\frac{\bar{\zeta}^{\beta}}{(1-\inner{w_j}{\zeta})^{|\beta|}}
\end{align*}
for all $z,\zeta\in\B_N$, where we define
\[\Bw(u)(z,\zeta) = (1-\inner{z}{\zeta})^{N+1+\gamma}\int_{\B_N}\frac{u(\xi)}{(1-\inner{z}{\xi})^{N+1+\gamma}(1-\inner{\xi}{\zeta})^{N+1+\gamma}}dV_{\gamma}(\xi)\]
Since the set 
\[\{1\}\cup\Big\{\frac{\bar{\zeta}^{\beta}}{(1-\inner{w_j}{\zeta})^{|\beta|}}:\quad 1\leq |\beta|\leq 2N+2\gamma,\quad 1\leq j\leq s\Big\}\]
is linearly independent, it follows that $h_1(z)$ and each $Q_{j,\beta}(\frac{z}{1-\inner{z}{w_j}})$ can be written as a linear combination of finitely many functions in the set
\[\Big\{\Bw(u)(\cdot,\zeta)-\bar{h}_2(\zeta): \quad\zeta\in \B_N\Big\}.\]
Note that for each $\zeta\in\B_N$, the function $\Bw(u)(z,\zeta)$ is the product of $(1-\inner{z}{\zeta})^{N+1+\gamma}$ with the weighted Bergman projection $\mathcal{P}_{\gamma}$ of $u(\xi)(1-\inner{\xi}{\zeta})^{-N-1-\gamma}$, which belongs to $L^{2N+2+2\gamma}_{\gamma}$ by the assumption about $u$. It is well known that $\mathcal{P}_{\gamma}$ maps $L^p_{\gamma}$ into $A^{p}_{\gamma}$ for $1<p<\infty$ (see \cite[Theorem~2.11]{Zhu2005}). Therefore, the function $\Bw(u)(\cdot,\zeta)-\bar{h}_2(\zeta)$ belongs to $L^{2N+2+2\gamma}_{\gamma}$. This implies that each $Q_{j,\beta}(\frac{z}{1-\inner{z}{w_j}})$ belongs to $L^{2N+2+2\gamma}_{\gamma}$. By Lemma \ref{L:L_2N_rational} below, for any $j$ with $w_j$ on the unit sphere, $Q_{j,\beta}$ must be constant, hence, identically zero since $Q_{j,\beta}$ vanishes at the origin. As a result, we may assume that $|w_j|<1$ for all $1\leq j\leq s$.

To complete the proof, we show that for $\omega\in\B_N$ and $1\leq j \leq N$, the rational function $\frac{z_j}{1-\inner{z}{\omega}}$ is a linear combination of $1$ and the components of $\varphi_{w}(z)$. The required representation then follows from \eqref{Eqn:B_u_formula}. 

For $z,\zeta\in\B_N$, \cite[Theorem~2.2.2]{Ru} provides the identity
\begin{align*}
    1-\inner{\varphi_{\omega}(z)}{\varphi_{\omega}(\zeta)} & = \frac{(1-|\omega|^2)(1-\inner{z}{\zeta})}{(1-\inner{z}{\omega})(1-\inner{\omega}{\zeta})},
\end{align*}
which is equivalent to
\[\frac{1-\inner{z}{\zeta}}{1-\inner{z}{\omega}} = \frac{1-\inner{\omega}{\zeta}}{1-|\omega|^2}\big(1-\inner{\varphi_{\omega}(z)}{\varphi_{\omega}(\zeta)}\big).\]
Setting $\zeta=0$
then $\zeta=e_j$ and subtracting the two quantities, we have
\[\frac{z_j}{1-\inner{z}{\omega}} = \frac{\omega_j}{1-|\omega|^2} + \frac{1}{1-|\omega|^2}\inner{\varphi_{\omega}(z)}{-\omega + (1-\omega_j)\varphi_{\omega}(e_j)}.\]
Note that the right hand-side is an affine function in $\varphi_{\omega}(z)$. As a consequence, for any multi-indexes $\alpha$ and $\beta$, the rational function \[\frac{z^{\alpha}\bar{z}^{\beta}}{(1-\inner{z}{\omega})^{|\alpha|}(1-\inner{\omega}{z})^{|\beta|}}\] is a polynomial  in $\varphi_{\omega}(z)$ and $\overline{\varphi_{\omega}(z)}$ of total degree $|\alpha|+|\beta|$.
\end{proof}

We now obtain a proof of statement (c) in Theorem \ref{main_Berezin}.
\begin{proposition}
\label{T:Bu_sum_L2N}
Suppose $u\in L^{1}_{\gamma}$ and $\Bw(u) = f_1\bar{g}_1+\cdots + f_d\bar{g}_d$, where $f_{\ell}, g_{\ell}\in \Hol$ and $f_{\ell}\in L^{2N+2+2\gamma}_{\gamma}$ for each $\ell$. Then there exist finitely many points $w_1,\ldots,w_s\in \B_{N}$, polynomials $Q_1,\ldots, Q_s$ in $\C[z,\bar{z}]$ with total degrees at most $2N+1+2\gamma$, and a pluriharmonic function $h$ such that
\[\Bw(u)(z) = h(z) + \sum_{j=1}^{s}Q_{j}\circ\varphi_{w_j}(z).\]
\end{proposition}

\begin{proof}
We know that there exist holomorphic functions $h_1, h_2$ on $\B_N$ and finitely many points $w_1,\ldots,w_s\in\overline{\B}_N$ and holomorphic polynomials $Q_{j,\beta}$ of degree at most $2N+1+2\gamma-|\beta|$ with $Q_{j,\beta}(0)=0$ such that
\begin{align*}
& h_1(z)+\bar{h}_2(z) + \sum_{\substack{1\leq j\leq s\\ \\
1\leq |\beta|\leq 2N+2\gamma}}Q_{j,\beta}\Big(\frac{z}{(1-\inner{z}{w_j})}\Big)\,\cdot\frac{\bar{z}^{\beta}}{(1-\inner{w_j}{z})^{|\beta|}}\\
&\qquad = \Bw(u)(z)\\
&\qquad = f_1(z)\,\overline{g_1(z)}+\cdots + f_d(z)\,\overline{g_d(z)}.
\end{align*}
Complexifying as in the proof of Theorem \ref{T:u_L2N} shows that each $Q_{j,\beta}\big(\frac{z}{1-\inner{z}{w_j}}\big)$ belongs to the linear span of
\[\{1\}\cup\Big\{\overline{g_1(\zeta)}\,f_1+\cdots+\overline{g_s(\zeta)}\,f_s:\ \zeta\in\B_N\Big\},\]
which is contained in $L^{2N+2+2\gamma}_{\gamma}$ by the hypothesis. The same argument as in the proof of Theorem \ref{T:u_L2N} may be used to finish the proof.
\end{proof}

\begin{lemma}
\label{L:L_2N_rational}
Let $Q$ be a polynomial in $\C[z_1,\ldots,z_N]$. If $Q\big(\frac{z}{1-\inner{z}{\omega}}\big)$ belongs to $L^{2N+2+2\gamma}_{\gamma}$ for some $\omega$ on the unit sphere, then $Q$ is a constant.
\end{lemma}

\begin{proof}
Since the case of a single complex variable may be regarded as a special case of two or more variables, we consider $N\geq 2$ throughout the proof. Without loss of generality, we may assume that $\omega=(0,\ldots,0,1)$. We write $z_{[N-1]}$ to denote $(z_1,\ldots,z_{N-1})\in\C^{N-1}$. Then $Q$ can be written as
\[Q(z) = \sum_{|\alpha|\geq 0}Q_{\alpha}(z_N)\,z_{[N-1]}^{\alpha},\] where the sum is finite over $\alpha\in\Z_{+}^{N-1}$ and each $Q_{\alpha}$ is a holomorphic polynomial in $z_N$. We have
\[F(z) = Q\Big(\frac{z}{1-\inner{z}{\omega}}\Big) = \sum_{|\alpha|\geq 0}Q_{\alpha}\Big(\frac{z_N}{1-z_N}\Big)\frac{z_{[N-1]}^{\alpha}}{(1-z_N)^{|\alpha|}}.\]
Since $F$ belongs to $L^{2N+2+2\gamma}_{\gamma}$, for each $\alpha$, the function
\begin{align*}
F_{\alpha}(z) & = Q_{\alpha}\Big(\frac{z_N}{1-z_N}\Big)\frac{z_{[N-1]}^{\alpha}}{(1-z_N)^{|\alpha|}}\\
& = \int_{[0,2\pi]^{N-1}}F(e^{i\theta_1}z_1,\ldots,e^{i\theta_{N-1}}z_{N-1},z_N)e^{-i(\alpha_1\theta_1+\cdots+\alpha_{N-1}\theta_{N-1})}\frac{d\theta_1}{2\pi}\cdots\frac{d\theta_{N-1}}{2\pi}
\end{align*}
must belong to $L^{2N+2+2\gamma}_{\gamma}$. However, for $|\alpha|\geq 1$, if $Q_{\alpha}$ is not identically zero, then for $z$ near $(0,\ldots,0,1)$, we see that $|F_{\alpha}(z)|$ dominates a nonzero constant multiple of $|g_{\alpha}(z)|$, where $g_{\alpha}(z) = \frac{z_{[N-1]}^{\alpha}}{(1-z_N)^{|\alpha|}}$.
Consider $z^{(\epsilon)} = (\epsilon,\ldots,\epsilon,1-N\epsilon^2)\rightarrow (0,\ldots,0,1)$ as $\epsilon\rightarrow 0$. A simple calculation reveals that $(1-|z^{(\epsilon)}|^2)^{1/2}|g(z^{(\epsilon)})|$ does \textit{not} converge to $0$ as $\epsilon\to 0$. By the remark after the proof of \cite[Theorem~2.1]{Zhu2005}, we conclude that $g_{\alpha}$ does not belong to $L^{2N+2+2\gamma}_{\gamma}$.

In addition, if $Q_{0}$ is not a constant, then for $z$ near $(0,\ldots,0,1)$, $|F_{0}(z)|=|Q_{0}(z)|$ dominates a nonzero constant multiple of $|\frac{1}{1-z_N}|$, which again does not belong to $L^{2N+2+2\gamma}_{\gamma}$. As a consequence, $Q_{\alpha}=0$ for all $|\alpha|\geq 1$ and $Q_{0}$ is a constant. Therefore, $Q$ is a constant as desired.
\end{proof}

Besides its important applications in the theory of Toeplitz operators as we shall see in the next section, Theorem \ref{main_Berezin} also helps answer open questions about $\mathcal{M}$-harmonic functions. In the early nineties, Ahern and Rudin \cite{AR} completely characterized holomorphic functions $f, g$ on the ball for which $f\bar{g}$ is $\mathcal{M}$-harmonic. Nearly a decade later, Zheng \cite{Zh} showed that for $f, g, h$ and $k$ belonging to the Hardy space $H^{2N}$, the function $f\bar{g}-h\bar{k}$ is $\mathcal{M}$-harmonic if and only if it is pluriharmonic. About ten years ago, making use of Ahern--Rudin's characterization, Choe et al. \cite[Lemma~4.5]{CKLJMAA2011} proved a single-product version of Zheng's result under a slightly weaker hypothesis. They only assumed that one of the factor belongs to $H^{2N}$. The problem of generalizing this and Zheng's result to finite sums of more than two products has been open since then, see \cite[Question~6.1]{CKLJMAA2011}. Using Theorem \ref{main_Berezin}, we obtain a far-reaching answer, in view of the fact that for all $-1< \delta < \gamma$ we have the continuous inclusions (\cite[Theorem ~5.13]{Bea89})
\[H^{2N} \subset A^{2N+2+2\delta}_{\delta} \subset A^{2N+2+2\gamma}_{\gamma}.\]

\begin{theorem}
\label{T:M_harmonic_open_question}
Suppose for each $1\leq j\leq s$, the functions $f_j, g_j$ are holomorphic on $\B_N$ and $f_j$ belongs to $A^{2N+2+2\gamma}_{\gamma}$. If $u = \sum_{j=1}^{s}f_j\bar{g}_j$ is an $L^{1}_{\gamma}$-eigenfunction of $\inL$, then $u$ must be pluriharmonic.
\end{theorem}

\begin{proof}
Using \cite[Theorem~4.2.4]{Ru}, it can be seen that $u$ is an eigenfunction of the Berezin transform $B_{\gamma}$, that is, there exists $\lambda\in\C$ such that
\[\Bw(u) = \lambda u = \sum_{j=1}^{s}\lambda f_j\bar{g}_j.\]
Since $u$ is clearly a $C^{2N+2+2\gamma}$-function and $f_j\in A^{2N+2+2\gamma}_{\gamma}$ for all $j$, Theorem \ref{main_Berezin} parts (a) and (c) hold, which implies that $u$ is pluriharmonic.
\end{proof}

\section{Brown--Halmos type results}
We first recall the following standard lemma characterizing when a function of the form $\sum_{j}\bar{g}_ju_j$ (with holomorphic $g_j, u_j$) is pluriharmonic. The one-dimensional version was already proved in \cite[Theorem~3.3]{CKL2008} but our proof here is much simpler. 

\begin{lemma}
\label{L:characterization_pluri}
Let $u_1,\ldots, u_s$ and $g_1,\ldots, g_s$ be holomorphic functions on $\B_N$. Then $\sum_{j=1}^{s}\bar{g}_ju_j$ is pluriharmonic if and only if
\[\sum_{j=1}^{s}\big(\overline{g}_j-\overline{g_j(0)}\big)\big(u_j-u_j(0)\big) = 0,\]
which is equivalent to
\[\sum_{j=1}^{s}\bar{g}_ju_j = \sum_{j=1}^{s}\Big(\bar{g}_ju_j(0) + \overline{g_j(0)}\,u_j-\overline{g_j(0)}u_j(0)\Big).\]
\end{lemma}

\begin{proof}
Without loss of generality, we may assume that $u_j(0)=g_j(0)=0$ for all $j$. Using power expansions, we have
\[\sum_{j=1}^{s}u_j(z)\bar{g}_j(z) = \sum_{|\alpha|\geq 1, |\beta|\geq 1}c_{\alpha,\beta}z^{\alpha}\bar{z}^{\beta},\] which is pluriharmonic if and only if it is identically zero.
\end{proof}

We are now ready to prove Theorem \ref{BH}, which is restated below for the reader's convenience.

\begin{theorem}\label{BH_restated}
Let $\phi_j, \psi_j$ for $1\leq j\leq n$ be bounded pluriharmonic functions 
and $h$ be a $C^{2N+2+2\gamma}$ bounded function on $\B_{N}$. Let $x_{\ell}, y_{\ell}\in \Aw$ for $1\leq\ell\leq r$. 
Write  $\phi_j=f_j+\bar{g}_j, \psi_j=u_j+\bar{v}_j$ where $f_j, g_j, u_j, v_j$ are holomorphic.
Then 
\begin{align}
\label{Eqn:finite_rank_rep}
\sum_{j=1}^nT_{\phi_j}T_{\psi_j}=T_h+\sum_{\ell=1}^{r}x_{\ell}\otimes y_{\ell}
\end{align}
if and only if $h- \sum_{j=1}^{n}\bar{g}_ju_j$ is pluriharmonic and
\begin{align}
\label{Eqn:formula_h}
\sum_{j=1}^{n} \phi_j\psi_j = h+(1-|z|^2)^{N+1+\gamma}\sum_{\ell=1}^{r}x_{\ell}\bar{y}_{\ell}.
\end{align}
\end{theorem}

\begin{proof}
For any functions $x,y\in \Aw$, we compute the Berezin transform
\[\Bw(x\otimes y)(z) = (1-|z|^2)^{N+1+\gamma}x(z)\,\overline{y(z)},\quad z\in\B_N.\]
Also, if $\phi=f+\bar{g}$ and $\psi=u+\bar{v}$ are bounded pluriharmonic, where $f, g, u, v$ are holomorphic functions (which might not be bounded but they all belong to $L^p_{\gamma}$ for all $p$), then it is well known that
\[\Bw(T_{\phi}T_{\psi}) = \phi\psi - \bar{g}u + \Bw(\bar{g}u).\]
Therefore,
\[\Bw\Big(\sum_{j=1}^{n}T_{\phi_j}T_{\psi_j}-T_h\Big) = \sum_{j=1}^{n}(\phi_j\psi_j - \bar{g}_ju_j) + \Bw\big(\sum_{j=1}^{n}\bar{g}_ju_j - h\big).\]
Using the linearity and injectivity of the Berezin transform, we conclude that \eqref{Eqn:finite_rank_rep} holds if and only if
\[\Bw\Big(\sum_{j=1}^{n}T_{\phi_j}T_{\psi_j} - T_h\Big) = \sum_{\ell=1}^{r}\Bw\big(x_{\ell}\otimes y_{\ell}\big),\]
which is equivalent to
\begin{align}
\label{Eqn:Ber_frk_sum_TOs}
\Bw\Big(\sum_{j=1}^{n}\bar{g}_ju_j - h\Big) = \sum_{j=1}^{n}(-\phi_j\psi_j+\bar{g}_ju_j)+(1-|z|^2)^{N+1+\gamma}\sum_{\ell=1}^{r}x_{\ell}\bar{y}_{\ell}.
\end{align}
We now show that this equation is equivalent to the two conditions stated in the theorem. Put $u=\sum_{j=1}^{n}\bar{g}_ju_j-h$. Suppose first that \eqref{Eqn:Ber_frk_sum_TOs} holds. Then the Berezin transform $\Bw(u)$ has finite rank because $\gamma$ is a integer. Since $u$ belongs to $C^{\infty}(\B_N)\cap L^{2N+2+2\gamma}_{\gamma}$, Theorem \ref{main_Berezin} implies that it is pluriharmonic on $\B_N$ and $\Bw(u)=u$. As a consequence, $h-\sum_{j=1}^{n}\bar{g}_ju_j$ is pluriharmonic and
\[\sum_{j=1}^{n}\bar{g}_ju_j-h = \sum_{j=1}^{n}(-\phi_j\psi_j+\bar{g}_ju_j)+(1-|z|^2)^{N+1+\gamma}\sum_{\ell=1}^{r}x_{\ell}\bar{y}_{\ell},\]
which gives \eqref{Eqn:formula_h}.

Conversely, if $u=\sum_{j=1}^{n}\bar{g}_ju_j-h$ is pluriharmonic and \eqref{Eqn:formula_h} holds, then using the fact that the Berezin transform fixes pluriharmonic functions, we conclude that \eqref{Eqn:Ber_frk_sum_TOs} holds, which implies \eqref{Eqn:finite_rank_rep} as desired.
\end{proof}

\begin{remark}
In this remark, we discuss a construction of functions that satisfy the two conditions in Theorem \ref{BH_restated}. As before, write  $\phi_j=f_j+\bar{g}_j, \psi_j=u_j+\bar{v}_j$, where $f_j, g_j, u_j, v_j$ are holomorphic and $f_j(0)=v_j(0)=0$. Assume that \eqref{Eqn:formula_h} holds, then
\begin{align*}
  h-\sum_{j=1}^{n}\bar{g}_ju_j = \sum_{j=1}^{n}(f_j\bar{v}_j+f_ju_j+\bar{g}_j\bar{v}_j) - (1-|z|^2)^{N+1+\gamma}\sum_{\ell=1}^{r}x_{\ell}\bar{y}_{\ell},
\end{align*}
which, by Lemma \ref{L:characterization_pluri}, is pluriharmonic if and only if 
\begin{align*}
    \sum_{j=1}^{n}f_j\bar{v}_j - (1-|z|^2)^{N+1+\gamma}\sum_{\ell=1}^{r}x_{\ell}\bar{y}_{\ell} & = \sum_{\ell=1}^{r}-x_{\ell}(0)\bar{y}_{\ell}-x_{\ell}\,\overline{y_{\ell}(0)}+x_{\ell}(0)\overline{y_{\ell}(0)}.
\end{align*}
The above identity is equivalent to
\begin{align}
\label{Eqn:pluriharmonic_identity_BH}
    \sum_{j=1}^{n}f_j(z)\bar{v}_j(z) & = \sum_{\ell=1}^{r}\bigg\{\big(x_{\ell}(z)-x_{\ell}(0)\big)\big(\bar{y}_{\ell}(z)-\overline{y_{\ell}(0)}\big)\\
    &\qquad\qquad + \sum_{1\leq|\alpha|\leq N+1+\gamma}(-1)^{|\alpha|}\binom{|\alpha|}{\alpha}\big(z^{\alpha}x_{\ell}(z)\big)\big(\overline{{z}^{\alpha}{y}_{\ell}(z)}\big)\bigg\}.\notag
\end{align}
Let $x_{\ell}, y_{\ell}$ ($1\leq\ell\leq r$) be any finite collection of bounded holomorphic functions. We can easily choose bounded holomorphic functions $f_j, v_j$ ($1\leq j\leq n$) for some $n$ such that $f_j(0)=v_j(0)=0$ and \eqref{Eqn:pluriharmonic_identity_BH} holds. For each $j$, choose arbitrary bounded holomorphic functions $g_j$ and $u_j$ and set $\phi_j=f+\bar{g}_j$ and $\psi_j = u_j+\bar{v}_j$. Put
\[h = \sum_{j=1}^{n}\phi_j\psi_j - (1-|z|^2)^{N+1+\gamma}\sum_{\ell=1}^{r}x_{\ell}\bar{y}_{\ell}.\]
We then have
\[\sum_{j=1}^{n}T_{\phi_j}T_{\psi_j} = T_{h} + \sum_{\ell=1}^{r}x_{\ell}\otimes y_{\ell}.\]
The problem becomes more delicate if one imposes a restriction on $n$. The paper \cite{DQZ} investigated the case $n=1$ in the setting of a single variable. It was shown that for bounded harmonic functions $\phi, \psi$, and smooth $h$, if $T_{\phi}T_{\psi}-T_h$ has rank one, then it must be zero. On the other hand, for any $r\geq 2$, examples were constructed so that $T_{\phi}T_{\psi}-T_h$ has rank exactly $r$. It would be interesting to generalize the results in \cite{DQZ} to the setting of several variables. 
\end{remark}

\begin{proof}[Proof of Corollary \ref{product}]
Write $\phi=f+\bar{g}, \psi=u+\bar{v}$ with holomorphic $f, g, u, v$ and $f(0)=v(0)=0$.

(a) By Theorem \ref{BH}, if $T_{\phi} T_{\psi}=T_h$, then $h=\phi\psi$
and  $h-\bar{g}u$ is pluriharmonic. It follows that $f\bar{v}=(h-\bar{g}u)-fu-\bar{g}\bar{v}$ is also pluriharmonic. Lemma \ref{L:characterization_pluri} implies that $f\bar{v}=0$ which forces either $f=0$ or $v=0$. Therefore, either $\bar{\phi}$ or $\psi$ must be holomorphic.

(b) Now suppose that $T_{\phi}T_{\psi}$  has a finite rank. Then there exist functions $x_{\ell}, y_{\ell}\in \Aw, 1\leq\ell\leq r$
so that $T_{\phi}T_{\psi}=\sum_{\ell} x_{\ell}\otimes y_{\ell}$. Using Theorem \ref{BH} with $h=0$, we obtain that 
 $\phi\psi=(1-|z|^2)^{N+1+\gamma}\sum_{\ell}x_{\ell}\bar{y}_{\ell}$,  and  $\bar{g}u$ is pluriharmonic, which implies either $g$ or $u$ is constant. Therefore, either $\phi$ or $\bar{\psi}$ is holomorphic. Taking operator adjoints if necessary, we may assume that $\phi$ is holomorphic. Assume further that $\phi$ is not identically zero. Then $T_{\phi}$ is injective. Since $T_{\phi}T_{\psi}$ has finite rank, it follows that $T_{\psi}$ must have finite rank, hence $\psi=0$ by the multivariable Luecking's Theorem.
\end{proof}

We now apply Theorem \ref{BH} to characterize when a sum of products of Hankel operators with pluriharmonic symbols has a finite rank. Recall that for a bounded symbol $\phi$, the Hankel operator $H_{\phi}: \Aw\rightarrow \Lw\ominus \Aw$ is defined as $H_{\phi}=(I-P_{\gamma})M_{\phi}|_{\Aw}$, where $M_{\phi}$ is the multiplication by $\phi$ and $P_{\gamma}$ is the weighted Bergman projection from $\Lw$ onto $\Aw$.
The crucial identity relating properties of Toeplitz and Hankel operators is given by
\[H_{\bar{\phi}}^{*}H_{\psi}=T_{\phi\psi}-T_{\phi}T_{\psi}.\]

\begin{proposition}
\label{P:finite_rank_sum_Hankel}
Let $\phi_j, \psi_j$ ($1\leq j\leq n$) be bounded pluriharmonic functions on $\B_N$. Then the following statements are equivalent:
\begin{enumerate}[(1)]
    \item $\sum_{j=1}^{n}H^{*}_{\bar{\phi}_j}H_{\psi_j} = 0$.
    \item $\sum_{j=1}^{n}H^{*}_{\bar{\phi}_j}H_{\psi_j} = T_{F}$ for some $F\in C^{2N+2+2\gamma}(\B_N)\cap L^{\infty}$.
    \item $\sum_{j=1}^{n}H^{*}_{\bar{\phi}_j}H_{\psi_j}$ has a finite rank.
    \item $\sum_{j=1}^{n}P(\phi_j)\cdot(\psi_j-P(\psi_j))$ is pluriharmonic.
\end{enumerate}
\end{proposition}

\begin{proof}
It is clear that (1) implies (2). Now assume that (2) holds. Then
\begin{align*}
    \sum_{j=1}^{n}T_{\phi_j}T_{\psi_j} = \sum_{j=1}^{n}(T_{\phi_j\psi_j}-H^{*}_{\bar{\phi}_j}H_{\psi_j}) = T_{h} - T_{F} = T_{h-F},
\end{align*}
where $h=\sum_{j=1}^{n}\phi_j\psi_j$. By Theorem \ref{BH}, we have
\[\sum_{j=1}^{n}\phi_j\psi_j = h - F,\] which implies $F=0$. Therefore, (3) (as well as (1)) follows.

Now assume that (3) holds, that is, the operator $T=\sum_{j=1}^{n}H^{*}_{\bar{\phi}_j}H_{\psi_j}$ has finite rank. The same argument as above gives $\sum_{j=1}^{n}T_{\phi_j}T_{\psi_j} = T_{h} - T$. By Theorem \ref{BH} again, the function
\[h-\sum_{j=1}^{n}(\phi_j-P(\phi_j))P(\psi_j)\] is pluriharmonic, which then implies (4).

Finally, assume that (4) holds. Setting $h=\sum_{j=1}^{n}\phi_j\psi_j$ and $x_{\ell}=y_{\ell}=0$, we see that both conditions in Theorem \ref{BH} are satisfied and so $\sum_{j=1}^{n}T_{\phi_j}T_{\psi_j} = T_{\sum_{j=1}^{n}\phi_j\psi_j}$, which gives (1). This completes the proof of the proposition.
\end{proof}

\begin{proof}[Proof of Corollary \ref{commutator}]
The sufficient direction is well known and not difficult to prove. To show the necessary direction, replacing $\phi$ by $\phi-\phi(0)$ and $\psi$ by $\psi-\psi(0)$ if necessary, we may assume that $\phi(0)=\psi(0)=0$. As before, write $\phi=f+\bar{g}$ and $\psi=u+\bar{v}$  with holomorphic $f, g, u,v$ satisfying $f(0)=g(0)=u(0)=v(0)=0$. We have
\[H_{\bar{\phi}}^{*}H_{\psi}-H_{\bar{\psi}}^{*}H_{{\phi}} = (T_{\phi\psi}-T_{\phi}T_{\psi})-(T_{\psi\phi}-T_{\psi}T_{\phi}) = -[T_{\phi}, T_{\psi}].\]
As a consequence, if $[T_{\phi}, T_{\psi}]$ has a finite rank, then so does $H_{\bar{\phi}}^{*}H_{\psi}-H_{\bar{\psi}}^{*}H_{{\phi}}$. Set $n=2$ and define
\[\phi_1 = \phi,\quad \phi_2 = -\psi,\quad \psi_1 = \psi \quad\text{and}\quad \psi_2 = \phi.\]
Then condition (3) in Proposition \ref{P:finite_rank_sum_Hankel} is satisfied. Therefore, (4) must hold, which means that
\[P(\phi)(\psi - P(\psi)) - P(\psi)(\phi-P(\phi)) = f\bar{v}-u\bar{g}\] is pluriharmonic. We may apply \cite[Theorem 5.6 and Lemma 6.8]{Zh} to complete the proof. Here, we provide a direct argument. Indeed, Lemma \ref{L:characterization_pluri} implies $f\bar{v}=u\bar{g}$ which, by complexifying, gives
\[f(z)\bar{v}(w) = u(z)\bar{g}(w)\quad\text{for all } z,w\in\B_N.\] 
If $u=v=0$, then $\psi=c\phi$ with $c=0$. If $u=0$ and $v$ is not identically zero, then $f=0$ so both $\phi$ and $\psi$ are anti-holomorphic. Similarly, if $v=0$ and $u$ is not identically zero, then $g=0$ so both $\phi$ and $\psi$ are holomorphic.
On the other hand, if neither of $u$ nor $v$ is identically zero, then there exists $z_0\in\B_N$ such that $u(z_0)v(z_0)\neq 0$ and it follows that $f=cu$ and $\bar{g}=c\bar{v}$, where
\[c = \frac{\bar{g}(z_0)}{\bar{v}(z_0)} = \frac{f(z_0)}{u(z_0)}.\]
Hence, $\phi-c\psi=0$. This completes the proof of the corollary.
\end{proof}

We end this section with another important application of Theorem \ref{BH}.

\begin{corollary}\label{BH-Rao}
Let $\phi_j, \psi_j\in L^{\infty}$ be pluriharmonic and let $h$ be $C^{2N+2+2\gamma}$-smooth and bounded. 
Write  $\phi_j=f_j+\bar{g}_j, \psi_j=u_j+\bar{v}_j$ where $f_j, g_j, u_j, v_j$ are holomorphic.
Then $\sum_{j=1}^nT_{\phi_j}T_{\psi_j}=T_h$ on $\Aw$ if and only if $h=\sum_{j=1}^{n}\phi_j\psi_j$ and
\[\sum_{j=1}^{n}\big(f_j-f_j(0)\big)\big(\bar{v}_{j}-\overline{v_{j}(0)}\big)=0.\]
\end{corollary}

\begin{remark}
\label{R:non_integer}
{It would be interesting to extend our results to general weighted Bergman spaces $A^2_{\gamma}$ with non-integer values $\gamma > -1$. However, there are two obstructions that we have not been able to resolve. First, let us recall that one of the main ingredients in our proof of Theorem \ref{main_Berezin} is the identity
\begin{align*}
|E_z|^2\Big(\frac{(1-|z|^2)^{N+1+\gamma}}{|1-\inner{z}{\xi}|^{2(N+1+\gamma)}}\Big) = \DO_{\xi}\Big(\frac{1}{|1-\inner{z}{\xi}|^2}\Big),
\end{align*}
where
\[\DO = \frac{1}{[(N+\gamma)!]^2}\Big(|E+(N+\gamma)|^2-\Delta\Big)\cdots\big(|E|^2-\Delta\big).\]
We do not know what $\DO$ should be if $\gamma$ is not an integer.

Second, in the proof of Theorem \ref{BH_restated}, we require that $\gamma$ be an integer so that the right hand-side of equation \eqref{Eqn:Ber_frk_sum_TOs} can be written as a finite sum of the form $\sum_{j}f_j\bar{g}_j$ with holomorphic $f_j$ and $g_j$. When $\gamma$ is not an integer, we have an infinite series so $B_{\gamma}(u)$ does not have a finite rank.

It would also be of interest to consider our results in the setting of pluriharmonic Bergman spaces. However, the situation is much less promising. Indeed, on the pluriharmonic Bergman space associated with the weighted measure $\dVw$, the reproducing kernel is given by
\[H^{\gamma}(w,z) = K^{\gamma}(w,z)+K^{\gamma}(z,w) - 1 = 2{\rm Re}\Big\{\frac{1}{(1-\inner{w}{z})^{N+1+\gamma}}\Big\} - 1\]
As a consequence, the Berezin transform admits a more complex integral formula. Our strategy does not seem suitable even when $\gamma = 0$.}
\end{remark}

\section{Polynomials in the range of Berezin transform and applications}
In this section, we investigate the range of the unweighted Berezin transform (denoted by $B$ instead of $B_0$ in this section). Due to the complicated computations, we are unable to consider the weighted case at this time. 
We first describe all polynomials in the range of $B$. We then construct examples which show that the conclusion of Theorem \ref{BH} may fail for $N\geq 2$ if the smoothness assumption on $h$ is dropped. Lastly, we show that the product of two Toeplitz operators with polynomial symbols, under a certain additional condition on the degrees, is always equal to another Toeplitz operator with an integrable symbol. We shall use $A^p$ and $L^p$ to denote the unweighted spaces.

In the setting of a single variable, Ahern \cite{A} showed that if $p$ and $q$ are holomorphic polynomials such that the degree of $pq$ is at most $3$, then $p\bar{q}$ is the Berezin transform of an $L^1$-function.
The following theorem generalizes this result to several variables. Since calculations cannot be performed explicitly as in the single variable case, the proof here is considerably more complicated. 

\begin{theorem} 
\label{P:polynomial_in_ran_B}
Let $f$ be a polynomials in $z$ and $\bar{z}$. Then $f=B(u)$ for some $u\in L^1$ if and only if for any $1\leq j, \ell\leq N$, the derivative $\partial_{z_j}\bar{\partial}_{z_{\ell}}f$ has total degree at most $2N-1$.

As a consequence, if $w_1,\ldots,w_s$ belongs to $\B_N$ and $Q_1,\ldots, Q_s$ are polynomials in $\C[z,\bar{z}]$ with total degrees at most $2N+1$, then there exists a function $u\in L^1$ such that $B(u) = \sum_{j=1}^{s}Q_j\circ\varphi_{w_j}$.
\end{theorem}

\begin{proof}
Throughout the proof, we write $\ran(B)$ to denote the image of $L^1$ under the Berezin transform.
Suppose that $f=B(u)$ for some $u\in L^{1}$. By Theorem \ref{T:Bu_sum_L2N}, there exist a pluriharmonic function $h$ and a polynomial $Q\in\C[z,\bar{z}]$ of degree at most $2N+1$ such that $f = B(u) = h + Q$. It follows that for any $1\leq j,\ell\leq N$,
\[\partial_{z_j}\bar{\partial}_{z_{\ell}}f = \partial_{z_j}\bar{\partial}_{z_{\ell}}h + \partial_{z_j}\bar{\partial}_{z_{\ell}}Q = \partial_{z_j}\bar{\partial}_{z_{\ell}}Q, \] which has total degree at most $2N-1$.

Conversely, suppose that for any $1\leq j, \ell\leq N$, the derivative $\partial_{z_j}\bar{\partial}_{z_{\ell}}f$ has total degree at most $2N-1$. Since $f$ is a polynomial, there exists a pluriharmonic polynomial $h$ and complex coefficients $c_{\alpha,\beta}$ for $|\alpha|\geq 1, |\beta|\geq 1$ such that
\[f(z) = h(z) + \sum_{|\alpha|\geq 1,|\beta|\geq 1}c_{\alpha,\beta}z^{\alpha}\bar{z}^{\beta}.\]
The assumption implies that $c_{\alpha,\beta}=0$ whenever $|\alpha|+|\beta| > 2N+2$. Consequently, we may write $f=h+Q$, where $h$ is pluriharmonic and $Q$ has total degree at most $2N+1$. Thus, it remains to show that $Q$ belongs to $\ran(B)$. 
 
Let $\alpha, \beta$ be two multi-indexes and $\ell$ be a non-negative integer such that $|\alpha|+|\beta|+2\ell\leq 2N+1$. We shall show that the polynomial $\bar{z}^{\alpha}z^{\beta}(1-|z|^2)^{\ell}$ belongs to $\ran(B)$. Taking complex conjugates if necessary, we may assume that $|\beta|\leq|\alpha|$. 

Using \cite[Proposition~1.4.9]{Ru} and the rotation invariant of the surface measure on $\partial\B_N$, we see that for any integer $s\geq 1$ and for any $z\in\B_N$,
\[\int_{\partial\B_N}|\inner{z}{\zeta}|^{2s}d\sigma(\zeta) = \frac{\Gamma(N)\,\Gamma(s+1)}{\Gamma(N+s)}|z|^{2s}.\]
Replacing $s$ by $s+|\alpha|$ and applying $\frac{\partial^{\alpha}_{z}}{(s+|\alpha|)\cdots (s+1)}$ to both sides of the above identity gives
\begin{align*}
\int_{\partial\B_N}\bar{\zeta}^{\alpha}\inner{z}{\zeta}^{s}\inner{\zeta}{z}^{s+|\alpha|}d\sigma(\zeta) = \frac{\Gamma(N)\,\Gamma(s+|\alpha|+1)}{\Gamma(N+s+|\alpha|)}\bar{z}^{\alpha}|z|^{2s}.
\end{align*}
Applying $\frac{\Gamma(s+|\alpha|-|\beta|+1)}{\Gamma(s+|\alpha|+1)}\bar{\partial}_z^{\beta}$, we have
\begin{align*}
    \int_{\partial\B_N}\bar{\zeta}^{\alpha}\zeta^{\beta}\inner{z}{\zeta}^{s}\inner{\zeta}{z}^{s+|\alpha|-|\beta|}d\sigma(\zeta) = \frac{\Gamma(N)\,\Gamma(s+|\alpha|-|\beta|+1)}{\Gamma(N+s+|\alpha|)}\bar{\partial}_z^{\beta}\Big(\bar{z}^{\alpha}|z|^{2s}\Big).
\end{align*}
Now let $u\in L^1$ be of the form $u(z) = \bar{z}^{\alpha}z^{\beta}\varphi(|z|^2)$, where $\varphi$ is a function on $[0,1)$ to be defined later. Integration in polar coordinates (using $\xi = r\zeta$) together with the above identity gives
\begin{align*}
    &\int_{\B_N}\bar{\xi}^{\alpha}\xi^{\beta}\varphi(|\xi|^2)\inner{z}{\xi}^{s}\inner{\xi}{z}^{s+|\alpha|-|\beta|}dV(\xi)\\
    &\qquad\qquad = 2N\int_{0}^{1}r^{2N+2s+2|\alpha|-1}\varphi(r^2)dr\int_{\partial\B_N}\bar{\zeta}^{\alpha}\zeta^{\beta}\inner{z}{\zeta}^{s}\inner{\zeta}{z}^{s+|\alpha|-|\beta|}d\sigma(\zeta)\\
    &\qquad\qquad = \frac{\Gamma(N+1)\,\Gamma(s+|\alpha|-|\beta|+1)}{\Gamma(N+s+|\alpha|)}\Big(\int_{0}^{1}r^{N+s+|\alpha|-1}\varphi(r)dr\Big)\bar{\partial}_z^{\beta}\big(\bar{z}^{\alpha}|z|^{2s}\big)\\
    &\qquad\qquad = \frac{\Gamma(N+1)\,\Gamma(s+|\alpha|-|\beta|+1)}{\Gamma(N+s+|\alpha|)}\widehat{\varphi}(N+s+|\alpha|)\,\bar{\partial}_z^{\beta}\big(\bar{z}^{\alpha}|z|^{2s}\big),
\end{align*}
where $\widehat{\varphi}$ denotes the Mellin transform of $\varphi$ given by
\[\widehat{\varphi}(\zeta) = \int_{0}^{1}r^{\zeta-1}\varphi(r)dr.\]
It follows that
\begin{align}
    \label{Eqn:integral_u}
   & \frac{1}{\Gamma(N+1)\,\Gamma(s+|\alpha|-|\beta|+1)}\int_{\B_N}u(\xi)\inner{z}{\xi}^{s}\inner{\xi}{z}^{s+|\alpha|-|\beta|}dV(\xi)\notag\\
   &\qquad\qquad\qquad = \frac{1}{\Gamma(N+s+|\alpha|)}\widehat{\varphi}(N+s+|\alpha|)\,\bar{\partial}_z^{\beta}\big(\bar{z}^{\alpha}|z|^{2s}\big).
\end{align}
We now compute, for $z\in\B_N$,
\begin{align*}
    & \int_{\B_N}\frac{u(\xi)}{|1-\inner{z}{\xi}|^{2(N+1)}}dV(\xi)\\
    & \ = \sum_{s,t=0}^{\infty}\frac{\Gamma(N+1+s)}{\Gamma(N+1)\,\Gamma(s+1)}\cdot\frac{\Gamma(N+1+t)}{\Gamma(N+1)\,\Gamma(t+1)}\int_{\B_N}u(\xi)\inner{z}{\xi}^{s}\inner{\xi}{z}^{t}dV(\xi).
\end{align*}
Since the integral vanishes unless $t=s+|\alpha|-|\beta|$, we may rewrite the above summation as
\begin{align}
\label{Eqn:integral_u_xi}
    & \int_{\B_N}\frac{u(\xi)}{|1-\inner{z}{\xi}|^{2(N+1)}}dV(\xi)\notag\\
    & \ = \sum_{s=0}^{\infty}\frac{\Gamma(N+1+s)}{\Gamma(N+1)\,\Gamma(s+1)}\cdot\frac{\Gamma(N+1+s+|\alpha|-|\beta|)}{\Gamma(N+1)\,\Gamma(s+|\alpha|-|\beta|+1)}\times\notag\\
    &\qquad\qquad\qquad\qquad\qquad\qquad\times\int_{\B_N}u(\xi)\inner{z}{\xi}^{s}\inner{\xi}{\xi}^{s+|\alpha|-|\beta|}dV(\xi)\\
    & \ = \bar{\partial}_z^{\beta}\left\{\bar{z}^{\alpha}\cdot\sum_{s=0}^{\infty}\frac{\Gamma(N+1+s)\,\Gamma(N+1+s+|\alpha|-|\beta|)}{\Gamma(N+1)\,\Gamma(s+1)\,\Gamma(N+s+|\alpha|)}\widehat{\varphi}(N+s+|\alpha|)|z|^{2s}\right\}.\notag
\end{align}
The last identity follows from formula \eqref{Eqn:integral_u}. To simplify the notation we now set $M=N+1-|\beta|-\ell$. 
 Since $|\alpha|+|\beta|+2\ell\leq 2N+1$ and $|\beta|\leq|\alpha|$, we have $1\leq M\leq N+1-|\beta|$. 
Let us choose $\varphi$ such that
\begin{align}
\label{Eqn:formula_Mellin_varphi}
\widehat{\varphi}(\zeta) 
& = \frac{\Gamma(N+1)}{\Gamma(N+1-|\ell|)}\cdot\frac{\Gamma(\zeta)\,\Gamma(\zeta+1-|\alpha|-|\beta|-\ell)}{\Gamma(\zeta+1-|\alpha|)\,\Gamma(\zeta+1-|\beta|)}\\
& = \frac{\Gamma(N+1)}{\Gamma(M+|\beta|)}\cdot\frac{\Gamma(\zeta)\,\Gamma(\zeta+M-N-|\alpha|)}{\Gamma(\zeta-|\alpha|+1)\,\Gamma(\zeta-|\beta|+1)}. \notag
\end{align}
The existence of such a function $\varphi$ will be established below. Since for all integers $s\geq 0$, 
\begin{align*}
    \widehat{\varphi}(N+s+|\alpha|) = \frac{\Gamma(N+1)}{\Gamma(M+|\beta|)}\cdot\frac{\Gamma(N+s+|\alpha|)\,\Gamma(M+s)}{\Gamma(N+s+1)\,\Gamma(N+s+|\alpha|-|\beta|+1)},
\end{align*}
formula \eqref{Eqn:integral_u_xi} simplifies to
\begin{align*}
    \int_{\B_N}\frac{u(\xi)}{|1-\inner{z}{\xi}|^{2(N+1)}}dV(\xi) & = \frac{\Gamma(M)}{\Gamma(M+|\beta|)}\notag \bar{\partial}_z^{\beta}\left\{\bar{z}^{\alpha}\cdot\sum_{s=0}^{\infty}\frac{\Gamma(M+s)}{\Gamma(M)\,\Gamma(s+1)}|z|^{2s}\right\}\\
    & = \frac{\Gamma(M)}{\Gamma(M+|\beta|)} \bar{\partial}_z^{\beta}\left\{\bar{z}^{\alpha}\,(1-|z|^2)^{-M}\right\}.
 \end{align*}
It follows that
\begin{align}
\label{Eqn:Berezin_formula_derivative}
    B(u)(z)
    & = \frac{\Gamma(M)}{\Gamma(M+|\beta|)}(1-|z|^2)^{N+1}\cdot \bar{\partial}_z^{\beta}\left\{\bar{z}^{\alpha}\,(1-|z|^2)^{-M}\right\}.
\end{align}
We now explain the existence of $\varphi$ and show that the corresponding function $u$ belongs to $L^1$. First, note that if $|\alpha|=0$, then $|\beta|=0$ as well since we assumed that $|\beta|\leq|\alpha|$ and so in this case, formula \eqref{Eqn:formula_Mellin_varphi} becomes
\begin{align*}
    \widehat{\varphi}(\zeta) 
    & = \frac{\Gamma(N+1)}{\Gamma(N+1-\ell)}\cdot\frac{\Gamma(\zeta)\,\Gamma(\zeta+1-\ell)}{\Gamma(\zeta+1)\,\Gamma(\zeta+1)}\\
    & = \begin{cases}
    \frac{1}{\zeta} & \text{ if } \ell=0,\\
    & \\
    \frac{\Gamma(N+1)}{\Gamma(N+1-\ell)}\cdot\frac{1}{\zeta}\cdot\frac{1}{(\zeta+1-\ell)\cdots \zeta} & \text{ if } \ell\geq 1.
    \end{cases}
\end{align*}
In the first case, $\varphi=1$. In the second case, $\varphi$ is a linear combination of $\log(r)$ and $r^{-1},\ldots, r^{1-\ell}$. 

Now assume $|\alpha|\geq 1$. Then second factor on the right hand-side of \eqref{Eqn:formula_Mellin_varphi} reduces to a proper rational function of the form
\[\begin{cases}
\frac{1}{(\zeta-|\beta|-\ell)\cdots (\zeta-|\beta|)} & \quad\text{ if } |\alpha|=1,\\
& \\
\frac{(\zeta+1-|\alpha|)\cdots (\zeta-1)}{(\zeta+1-|\alpha|-|\beta|-\ell)\cdots (\zeta-|\beta|)} & \quad\text{ if } |\alpha|\geq 2,
\end{cases}\]
whose numerator has degree $|\alpha|-1$ and whose denominator has degree $|\alpha|+\ell > |\alpha|-1$. Therefore, $\varphi(r)$ exists and it is a linear combination of $r^{1-|\alpha|-|\beta|-\ell}, \ldots, r^{-|\beta|}$. 
In all cases, we have
\[\varphi(r) = O(r^{1-|\alpha|-|\beta|-\ell})\quad\text{ as } r\to 0^{+},\]
which implies that for any $\zeta\in\partial\B_N$,
\[u(r\zeta) = r^{|\alpha|+|\beta|}\varphi(r^2)\bar{\zeta}^{\alpha}\zeta^{\beta} = O(r^{2-|\alpha|-|\beta|-2\ell}).\]
Since $(2N-1)+2-|\alpha|-|\beta|-2\ell=2N+1-|\alpha|-|\beta|-2\ell\geq 0$, using integration by polar coordinates, we conclude that $u\in L^1$.

Choosing $|\beta|=0$ in \eqref{Eqn:Berezin_formula_derivative} shows that $\bar{z}^{\alpha}(1-|z|^2)^{\ell}$ belongs to $\ran(B)$ whenever $|\alpha|+2\ell\leq 2N+1$. It then follows that $\bar{z}^{\alpha}|z|^{2s}$ (and hence $z^{\alpha}|z|^{2s}$, after taking complex conjugates) belongs to $\ran(B)$ whenever $|\alpha|+2s\leq 2N+1$. 

Generally, whenever $|\alpha|+|\beta|+2\ell\leq 2N+1$, we may use \eqref{Eqn:Berezin_formula_derivative} to conclude that $\ran(B)$ contains the function
\begin{align*}
     & \frac{\Gamma(M)}{\Gamma(M+|\beta|)}(1-|z|^2)^{N+1}\cdot \bar{\partial}_z^{\beta}\left\{\bar{z}^{\alpha}\,(1-|z|^2)^{-M}\right\}\\
    &\qquad\qquad = \bar{z}^{\alpha}z^{\beta}(1-|z|^2)^{\ell} + \sum_{\substack{\mu+\nu=\beta\\ |\mu|\geq 1}}c_{\mu,\nu}\,\bar{z}^{\alpha-\mu}\,z^{\nu}\,(1-|z|^2)^{N+1-M-|\nu|}\\
    &\qquad\qquad = \bar{z}^{\alpha}z^{\beta}(1-|z|^2)^{\ell} + \sum_{\substack{\mu+\nu=\beta\\ |\mu|\geq 1}}c_{\mu,\nu}\,\bar{z}^{\alpha-\mu}\,z^{\nu}\,(1-|z|^2)^{|\beta|+\ell-|\nu|},
\end{align*}
where $c_{\mu,\nu}$'s are constants. Note that each term in the summation has total degree at most $|\alpha|+|\beta|+2\ell\leq 2N+1$ and the degree in $z$ is $|\nu| < |\beta|$. As a consequence, an induction in $|\beta|$ shows that $\bar{z}^{\alpha}z^{\beta}(1-|z|^2)^{\ell}$ belongs to $\ran(B)$ whenever $|\alpha|+|\beta|+2\ell\leq 2N+1$. Letting $\ell=0$, we conclude that $\bar{z}^{\alpha}z^{\beta}\in\ran(B)$ whenever $|\alpha|+|\beta|\leq 2N+1$. 

For each $1\leq j\leq s$, we showed above the existence of a function $u_j\in L^1$ such that $B(u_j)=Q_j$. Using the commutativity of the Berezin transform and automorphisms of the unit ball (see \cite[Proposition~2.3]{AFR}, for example), we have
\[B(u_j\circ\varphi_{w_j}) = B(u_j)\circ\varphi_{w_j} = Q_j\circ\varphi_{w_j}.\]
It then follows that $B(\sum_{j=1}^{s}u_j\circ\varphi_{w_j}) = \sum_{j=1}
^{s}Q_j\circ\varphi_{w_j}$ as required.
\end{proof}

\begin{remark}
Using \eqref{Eqn:Berezin_formula_derivative} in the case $\ell=1$, $|\beta|=0$ and $|\alpha|\geq 1$ (hence $M=N$) shows that for $u(z)=\bar{z}^{\alpha}\varphi(|z|^2)$ with $\widehat{\varphi}$ given by \eqref{Eqn:formula_Mellin_varphi}, we obtain
\begin{align*}
   B\Big(\bar{z}^{\alpha}|z|^{-2|\alpha|}-\bar{z}^{\alpha}\Big)  = \frac{|\alpha|}{N}\bar{z}^{\alpha}(1-|z|^2),
\end{align*}
which implies
\[B\Big(\frac{|\alpha|+N}{N}\bar{z}^{\alpha}-\bar{z}^{\alpha}|z|^{-2|\alpha|}\Big) = \frac{|\alpha|}{N}\bar{z}^{\alpha}|z|^2.\]
This identity is valid for all $1\leq |\alpha|<2N$. For $N=1$ and $\alpha=1$, this is the formula given in \cite[Lemma~1]{A}. 
\end{remark}

As it is well known in the literature, properties of the Berezin transform have consequences in the theory of Toeplitz operators. We discuss here a few examples.
\begin{remark}
Setting $\alpha=\beta=(1,0,\ldots,0)$ and $\ell=0$ (hence, $M=N$) in \eqref{Eqn:Berezin_formula_derivative} gives
\begin{align}
\label{Eqn:Bu_bounded}
B(u)(z) & =\frac{1}{N}(1-|z|^2)^{N+1}\cdot\bar{\partial}_{z_1}\{\bar{z}_1(1-|z^2|)^{-N}\}\notag\\
& = \frac{1}{N}(1-|z|^2)^{N+1}\cdot\big\{(1-|z|^2)^{-N}+N|z_1|^2(1-|z|^2)^{-N-1}\big\}\notag\\
& = \frac{1}{N}(1-|z|^2) + |z_1|^2.
\end{align}
Here, $u(z)=|z_1|^2\varphi(|z|^2)$ with
\begin{align*}
    \widehat{\varphi}(\zeta) & = \frac{\Gamma(\zeta)\,\Gamma(\zeta-1)}{\Gamma(\zeta)\,\Gamma(\zeta)} = \frac{1}{\zeta-1} = \widehat{r^{-1}}(\zeta).
\end{align*}
It follows that $u(z)=\frac{|z_1|^2}{|z|^2}$ for $z\in\B_N\backslash\{0\}$, which is bounded and is \textit{not} pluriharmonic if $N\geq 2$. 

Let us rewrite the identity in \eqref{Eqn:Bu_bounded} in the form
\[(N-1)|z_1|^2 - \sum_{j=2}^{N}|z_j|^2 = B(h)\]
with $h(z) = -1+\frac{N|z_1|^2}{|z|^2}$. It follows that
\[(N-1)T_{z_1}T_{\bar{z}_1} - \sum_{j=2}^{N}T_{z_j}T_{\bar{z}_j} = T_{h},\]
where $h$ is a bounded function and $h(z)\neq (N-1)|z_1|^2-\sum_{j=2}^{N}|z_j|^2$. It is important to note that this phenomenon cannot occur for $N=1$ due to Corollary \ref{BH-Rao}.
\end{remark}

We end the paper by showing that the product of two Toeplitz operators with polynomial symbols, under a certain condition on the degrees, is again a Toeplitz operator. However, we note that the symbol of the resulting Toeplitz operator is not always a polynomial.

\begin{proposition}
\label{T:product_monomial_Toeplitz}
Let $\alpha$ and $\beta$ be two multi-indexes such that $|\alpha|\geq 1$ and $|\beta|\geq 1$. Then $T_{z^{\beta}}T_{\bar{z}^{\alpha}} = T_{u}$ for some $u\in L^1$ if and only if $|\alpha|+|\beta|\leq 2N+1$.

As a consequence, if $f$ and $g$ are polynomials in $z$ and $\bar{z}$ such that the sum of the degree of $f$ in $z$ and the degree of $g$ in $\bar{z}$ is at most $2N+1$, then there exists $h\in L^1$ such that $T_{f}T_{g} =T_h$.
\end{proposition}

\begin{proof}
Suppose $T_{z^{\beta}}T_{\bar{z}^{\alpha}} = T_{u}$ for some $u\in L^1$. Taking Berezin transforms gives
\[B(u) = B(T_{u}) = B(T_{z^{\beta}}T_{\bar{z}^{\alpha}}) = \bar{z}^{\alpha}z^{\beta}.\]
Write $\alpha=(\alpha_1,\ldots,\alpha_N)$ and $\beta=(\beta_1,\ldots,\beta_N)$. Since $|\alpha|\geq 1$ and $|\beta|\geq 1$, there exist $j, \ell$ such that $\beta_{j}\neq 0$ and $\alpha_{\ell}\neq 0$. It follows that the total degree of $\partial_{z_{j}}\bar{\partial}_{z_{\ell}}(\bar{z}^{\alpha}z^{\beta})$ is exactly $|\alpha|+|\beta|-2$. Proposition \ref{P:polynomial_in_ran_B} implies that $|\alpha|+|\beta|-2\leq 2N-1$, which gives $|\alpha|+|\beta|\leq 2N+1$.

Conversely, if $|\alpha|+|\beta|\leq 2N+1$, then by Theorem \ref{P:polynomial_in_ran_B}, there exists
a function $u\in L^1$ such that $\bar{z}^{\alpha}z^{\beta}=B(u)$. 
This implies that $B(T_{z^{\beta}}T_{\bar{z}^{\alpha}}) = B(u)$, which gives $T_{z^{\beta}}T_{\bar{z}^{\alpha}}=T_{u}$. 
For any holomorphic polynomials $p$ and $q$, using the well-known properties of Toeplitz operators, we have
\[T_{\bar{p}(z)z^{\beta}}T_{q(z)\bar{z}^{\alpha}} = T_{\bar{p}}\big(T_{z^{\beta}}T_{\bar{z}^{\alpha}}\big)T_{q}  = T_{\bar{p}}T_{u}T_{q} = T_{\bar{p}uq}.\]
The last statement of the proposition now follows.
\end{proof}

\medskip

\textit{Acknowledgements}. We would like to thank the referee for a thorough reading and many constructive comments that improved the presentation of the paper. In particular, it was the referee's request that motivated us to extend our prior results to the weighted Bergman spaces.

\bibliography{ref}
\bibliographystyle{alpha}
\end{document}